\documentclass[aps,prl,twocolumn,floatfix,nofootinbib]{revtex4-1}
\usepackage{graphics}
\usepackage{bm}
\usepackage{subfigure}
\usepackage{graphicx}
\usepackage{amsthm}
\usepackage{notes2bib}
\usepackage{amsmath}
\usepackage{amssymb}
\usepackage{url}
\usepackage{enumerate}
\usepackage{epstopdf}
\usepackage{color}
\usepackage{mathtools}          
\usepackage{setspace}

\usepackage[usenames,dvipsnames,svgnames,table]{xcolor}
\usepackage{color}

\newcommand{\mb}{\mathbf}

\usepackage{lineno}
\usepackage{comment}
\begin{document}
\title{
Dynamical system analysis of a data-driven model constructed by reservoir computing
}
\author{Miki U. Kobayashi}
\affiliation{Faculty of Economics, Rissho University, Tokyo 141-8602, Japan}
\author{Kengo Nakai}
\affiliation{Faculty of Marine Technology,
   Tokyo University of Marine Science and Technology,  Tokyo 135-8533, Japan}
\author{Yoshitaka Saiki}
\affiliation{Graduate School of Business Administration, Hitotsubashi University, Tokyo 186-8601, Japan}
\author{Natsuki Tsutsumi}
\affiliation{Faculty of Commerce and Management, Hitotsubashi University, Tokyo 186-8601, Japan}
\begin{abstract}
This study evaluates data-driven models from a dynamical system perspective, such as unstable fixed points, periodic orbits, chaotic saddle, Lyapunov exponents, manifold structures, and statistical values. 
We find that these dynamical characteristics can be reconstructed much more precisely by a data-driven model 
than by computing directly from training data.
With this idea, we predict the laminar lasting time distribution of a particular macroscopic variable of chaotic fluid flow, 
which cannot be calculated from a direct numerical simulation of the Navier--Stokes equation because of its high computational cost. 
\end{abstract}
\date{\today}
\maketitle

\section{I. Introduction.}
Reservoir computing, a brain-inspired machine-learning technique that employs a data-driven dynamical system,    
is effective in predicting time series and frequency spectra 
in chaotic behaviors, including fluid flow and 
global atmospheric dynamics~\cite{Verstraeten_2007,
Zhixin_2017, Pathak_2017,
Pathak_2018,Antonik_2018,nakai_2018,arcomano_2020,pandey_2020,huang_2020,kong_2021}.
\textcite{Pathak_2017} examined the Lorenz system and the Kuramoto--Sivashinsky system and reported that the data-driven model obtained from reservoir computing 
could generate an arbitrarily long time series that mimics the dynamics of the original systems.\\ 
\indent The extent to which a data-driven model using reservoir computing can capture the dynamical properties of original systems should be determined. 
\textcite{Lu_2018} reported that 
a data-driven model has an attractor similar to that of the original system under an appropriate choice of parameters.
\textcite{nakai_2020} confirmed that a single data-driven model could infer the time series of chaotic fluid flow from various initial conditions.
\textcite{Zhu_2019} identified 
some unstable periodic orbits of a  data-driven model through delayed feedback control.
They suggested that a data-driven model could reconstruct the attractor of the original dynamical system.\\
\indent This paper clarifies that a data-driven model using reservoir computing has richer information than that obtained from a training data, especially from dynamical system point of view, 
suggesting that dynamical properties of 
the original unknown dynamical system can be estimated 
by reservoir computing from a relatively short time series.
Besides the invariant sets, such as fixed points and periodic orbits, 
the dynamical properties, such as Lyapunov exponents and  manifold structures between stable and unstable manifolds, can be reconstructed by the data-driven model through reservoir computing,
 even if the system does not have structural stability.\\
\indent We mainly deal with the 
Lorenz system~\cite{l1963}:
\begin{equation}
\frac{dx}{dt}=10 (y-x),\  \frac{dy}{dt}=rx-y-xz,\ \frac{dz}{dt}=xy-\frac{8}{3}z,\label{eq:lorenz}
\end{equation}
and will be denoted as the actual Lorenz system in this paper.
A data-driven model is constructed from a short time training data created from (\ref{eq:lorenz}), the method of which is explained later.
Two different parameter values of $r$ are considered. One of the parameters ($r=28$) 
has hyperbolic dynamics, 
whereas the other ($r=60$) generates dynamics 
with tangencies between stable and unstable manifolds~\cite{saiki_2010}.
The latter property is one of the two primary sources for the breaking structural stability~\cite{bonatti_2005}, which often appears in the real-world physical phenomena.
We also deal with the R\"ossler system~\cite{roessler_1976}:
\begin{equation}
\frac{dx}{dt}=- y - z,\  \frac{dy}{dt}=x + 0.2y,\ \frac{dz}{dt}=0.2 + ( x - 5.7 )  z,\label{eq:rossler}
\end{equation}
in order to confirm that the similar properties hold. 
As an application of the obtained knowledge, this study examines high-dimensional chaotic fluid flow to determine if 
the laminar lasting time 
distribution can be predicted using the data-driven model constructed from short training time-series data. \\
\indent 
After introducing the method of reservoir computing in Section II, we investigate the dynamical system properties of the data-driven model obtained from the reservoir computing 
for the Lorenz system in Section III and the R\"ossler system in Section IV. Applying the obtained implications, in Section V, we estimate the state-lasting time  distribution. We conclude our remarks in Section VI.
\section{II. Reservoir computing.}
A reservoir is a recurrent neural network whose internal parameters are not adjusted to fit the data in the training process~\cite{Jaeger_2001,Jaeger_2004}.
The reservoir can be trained by feeding it an input time series and fitting a linear function of the reservoir state vector 
to the desired output time series. 
We do not use a physical knowledge in constructing a model. 
The data-driven model using reservoir computing we study is the following:
\begin{equation}
\begin{cases}
    \mb{u}(t)=\mb{W}^*_{\text{out}}\mb{r}(t), \\
	\mb{r}(t+\Delta t)=(1-\alpha)\mb{r}(t)+\alpha \tanh(\mb{A}\mb{r}(t)+\mb{W}_{\text{in}}\mb{u}(t)
	),
	\label{eq:reservoir}
	\end{cases}
	\end{equation}
where 
$\mb{u}(t) \in \mathbb{R}^M$ is a vector-valued variable, 
the component of which is  
denoted as an output variable; 
$\mb{r}(t) \in \mathbb{R}^N~(N \gg M)$ is a reservoir state vector; 
$\mb{A} \in \mathbb{R}^{N\times N}$, 
$\mb{W}_{\text{in}} \in \mathbb{R}^{N\times M}$, 
and $\mb{W}^*_{\text{out}} \in \mathbb{R}^{M\times N}$
are matrices;
$\alpha$ ($0<\alpha\le 1$) is a coefficient; 
$\Delta t$ is a time step. 
We define $\tanh(\mb{q})=(\tanh(q_1), \tanh(q_2),\ldots,\tanh(q_N))^{\text{T}},$
for a vector $\mb{q} = (q_1,q_2,\ldots,q_N)^{\text{T}}$, 
where $\text{T}$ represents the transpose of a vector.

We explain how to determine $\mb{W}^*_{\text{out}}$ in (\ref{eq:reservoir}). 
Time development of the reservoir state vector 
$\mb{r}(l \Delta t)$ are determined by 
\begin{equation}
\mb{r}(t+\Delta t)=(1-\alpha)\mb{r}(t)+\alpha \tanh(\mb{A}\mb{r}(t)+\mb{W}_{\text{in}}\mb{u}(t)
	),\label{eq:r-iteration}
\end{equation}
together with training time-series data $\{\mb{u}(l \Delta t)\} (-L_0\le l \le L)$, where $L_0$ is the transient time and  $L$ is the time length to determine $\mb{W}^*_{\text{out}}$. 
For given random matrices 
$\mb{A}$ and $\mb{W}_{\text{in}}$, 
we determine $\mb{W}_{\text{out}}$ so that
the following quadratic form takes the minimum: 
\begin{equation}
\displaystyle\sum^{L}_{l=0} \|\mb{W}_\text{out}\mb{r}(l\Delta t)-\mb{u}((l+1)\Delta t)\|^2
+\beta[Tr(\mb{W}_\text{out}\mb{W}^{\text{T}}_\text{out})],\label{eq:minimize}
\end{equation}
where $\|\mb{q}\|^2=\mb{q}^{\text{T}} \mb{q}$ for a vector $\mb{q}$. 
The minimizer is 
\begin{align}
	\mb{W}^*_\text{out}&=\delta\mb{U}\delta\mb{R}^{T}(\delta\mb{R}\delta\mb{R}^{T}+\beta\mb{I})^{-1}, \label{eq:wout-c1}
\end{align}
where 
$\mb{I}$ is the $N \times N$ identity matrix, $\delta\mb{R}$ (respectively, $\delta\mb{U}$) is the matrix
whose $l$-th column is $\mb{r}(l\Delta t)$ (respectively, $\mb{u}(l\Delta t)$). 
(see~\cite{Lukosevicius_2009}~P.140 and \cite{Tikhonov_1977} Chapter 1  for details).

Note that $\mb{A}$ is chosen to have a maximum eigenvalue $\rho$ $(|\rho|<1)$ 
in order for (\ref{eq:r-iteration}) to satisfy so called echo state property.
It is known that adding noise to the training time-series data can be useful in the construction of a data-driven model \cite{Zhixin_2017}. 
{For the computation of the data-driven model of the R\"ossler system, a small amplitude of noise is added.}
More details about the reservoir computing can be found elsewhere~\cite{Pathak_2017,nakai_2018}.\\

\begin{table}[htb]
\small
	\begin{tabular}{|l|l|r|r|}
	\hline
    \multicolumn{2}{|c|}{parameter}&$r=28$&$r=60$ \\ \hline
    ~$M$&dimension of input and output variables& \multicolumn{2}{c|}{6} \\  \hline
	~$N$&dimension of reservoir state vector & \multicolumn{2}{c|}{2000} \\  \hline
    ~$\Delta t$&time step for a model~(\ref{eq:reservoir}) & \multicolumn{2}{c|}{0.01} \\  \hline
    ~$\rho$&maximal eigenvalue of $\mb{A}$ &\multicolumn{2}{c|}{0.99} \\  \hline
    ~$\alpha$&nonlinearity degree in a model~(\ref{eq:reservoir})&0.3&0.4  \\ \hline
    ~$\beta$&regularization parameter &0.002&0.001 \\ \hline
    ~$\Delta \tau$&delay-time for input and output variables &0.11&0.07 \\ \hline
    \end{tabular}
    \caption{{\bf The list of parameters and their values used in the reservoir computing in each section.}
    We use
$\mb{u}(t) = (x(t),y(t),z(t),x(t-\Delta \tau),y(t-\Delta \tau),z(t-\Delta \tau))$ for the input variable, where $\Delta\tau$ is the delay time.
 	}
 	\label{tab:parameter}
\end{table}

\section{III. Lorenz system}
\indent In this section we evaluate a data-driven model~(\ref{eq:reservoir}) constructed using short training time series data from a dynamical system perspective.
The main focus is on the properties in the space of output variables (corresponding to $x$, $y$ and $z$ for the case of the Lorenz system),  
which compare them with those of the actual system.
The sets of parameter values used to construct the data-driven model are shown in Table \ref{tab:parameter}.
%
     \begin{figure}
		\subfigure[\hspace{-2mm}]{\includegraphics[width=0.435\columnwidth,height=0.45\columnwidth]{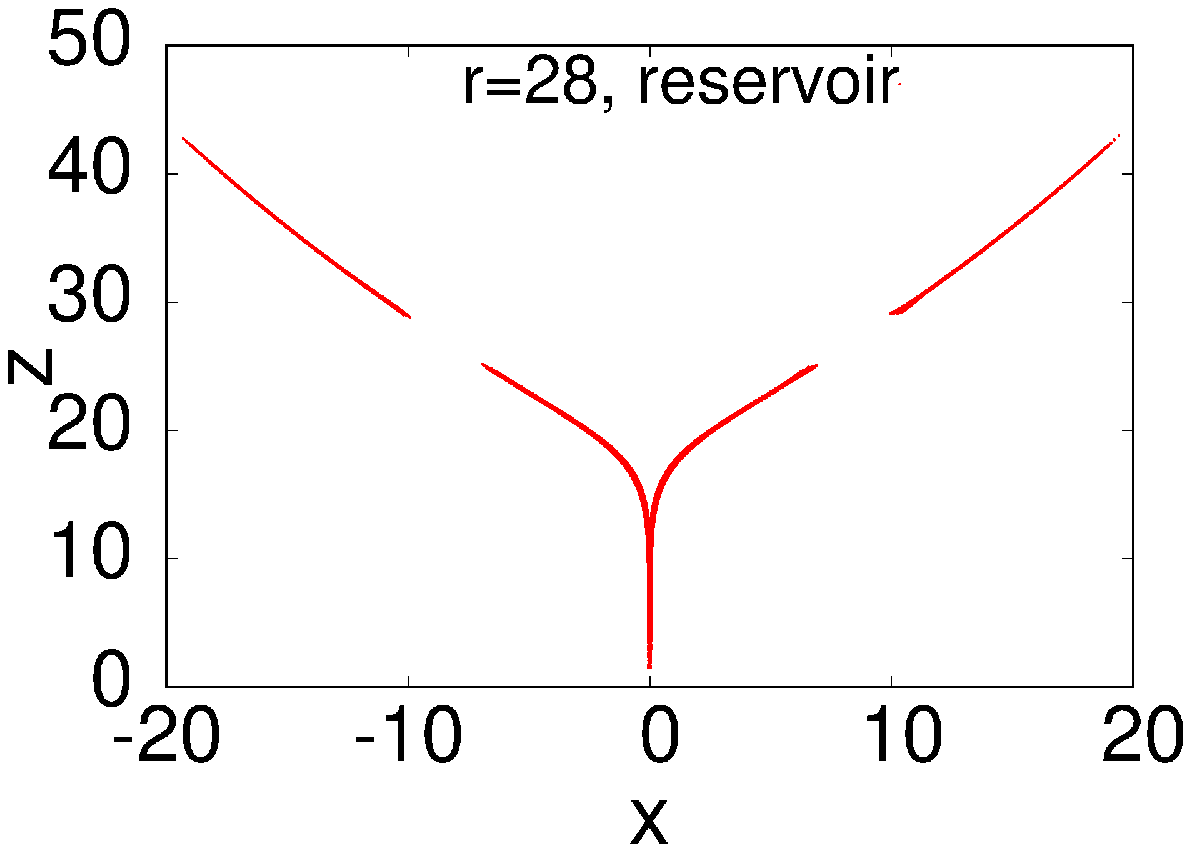}}%
		\subfigure[\hspace{-4mm}]{\includegraphics[width=0.455\columnwidth,height=0.45\columnwidth]{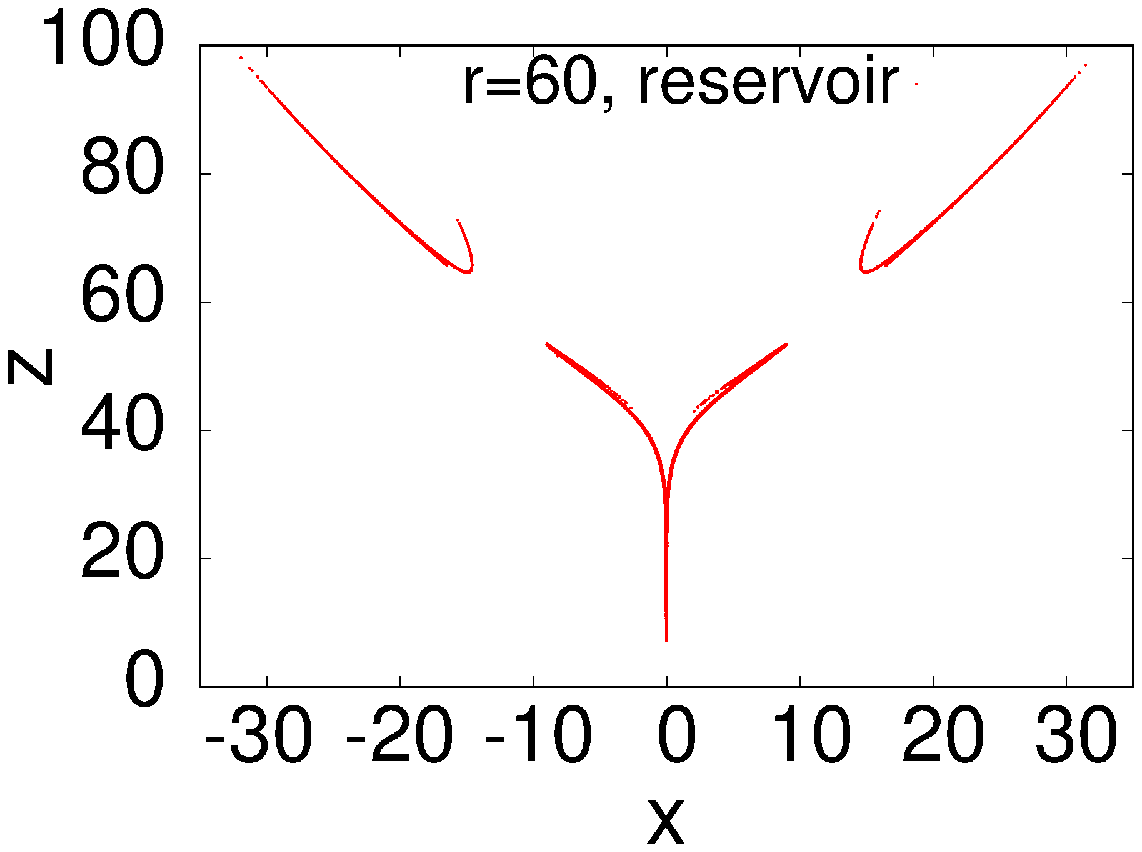}}\\%
		\subfigure[\hspace{-2mm}]{\includegraphics[width=0.435\columnwidth,height=0.45\columnwidth]{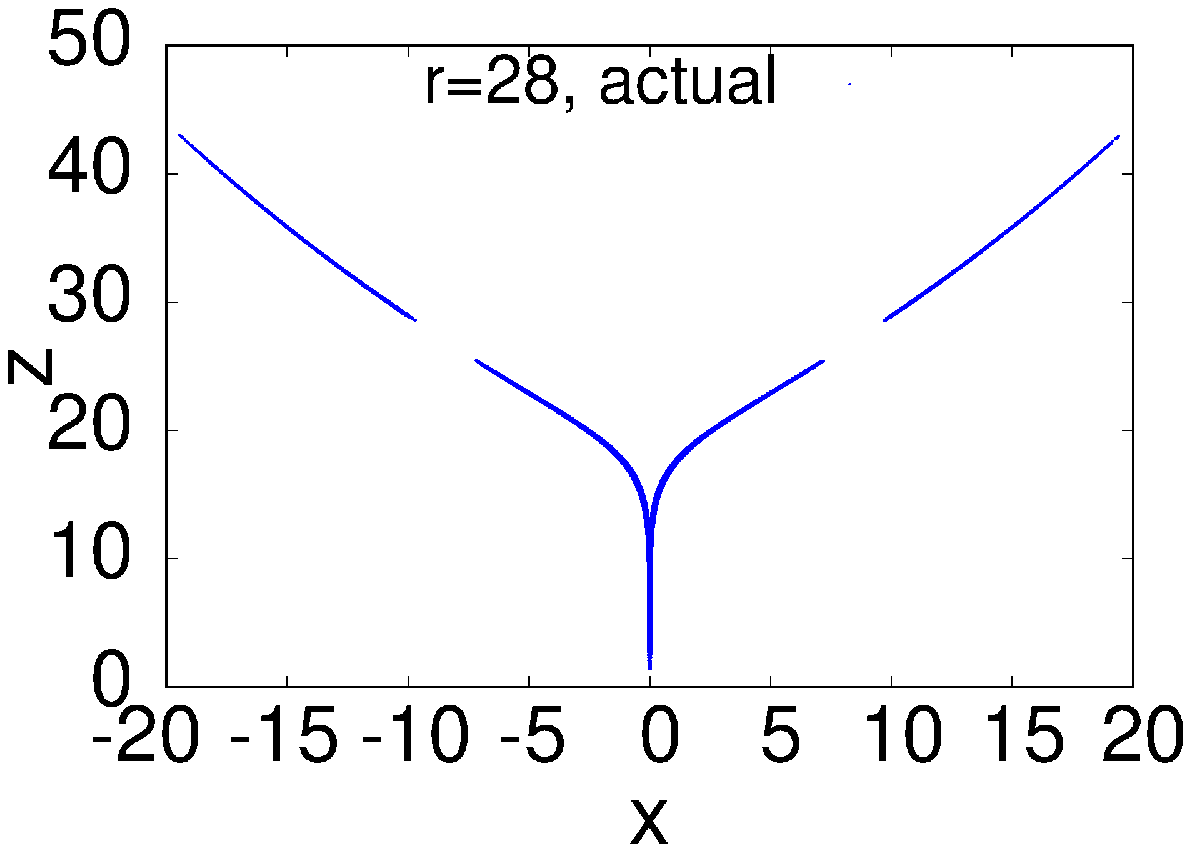}}%
		\subfigure[\hspace{-4mm}]{\includegraphics[width=0.455\columnwidth,height=0.45\columnwidth]{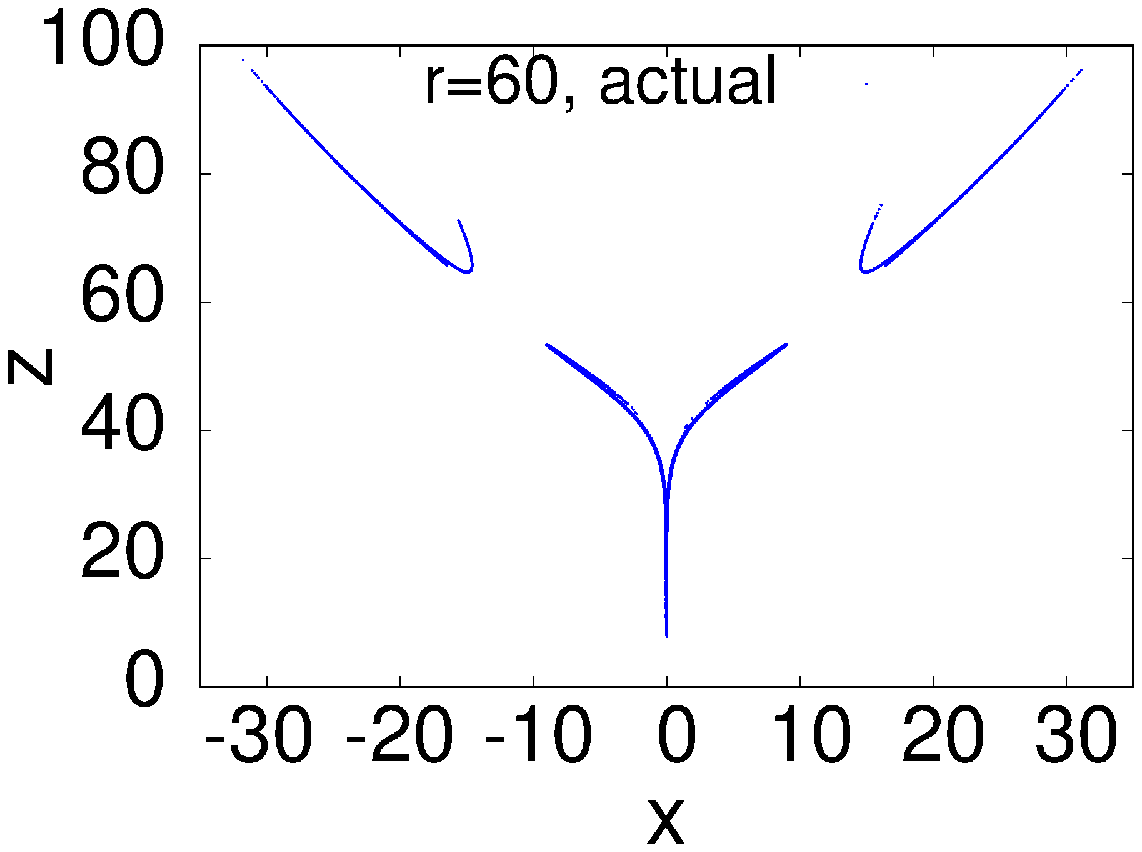}}\\%
		\subfigure[\hspace{-2mm}]{\includegraphics[width=0.435\columnwidth,height=0.45\columnwidth]{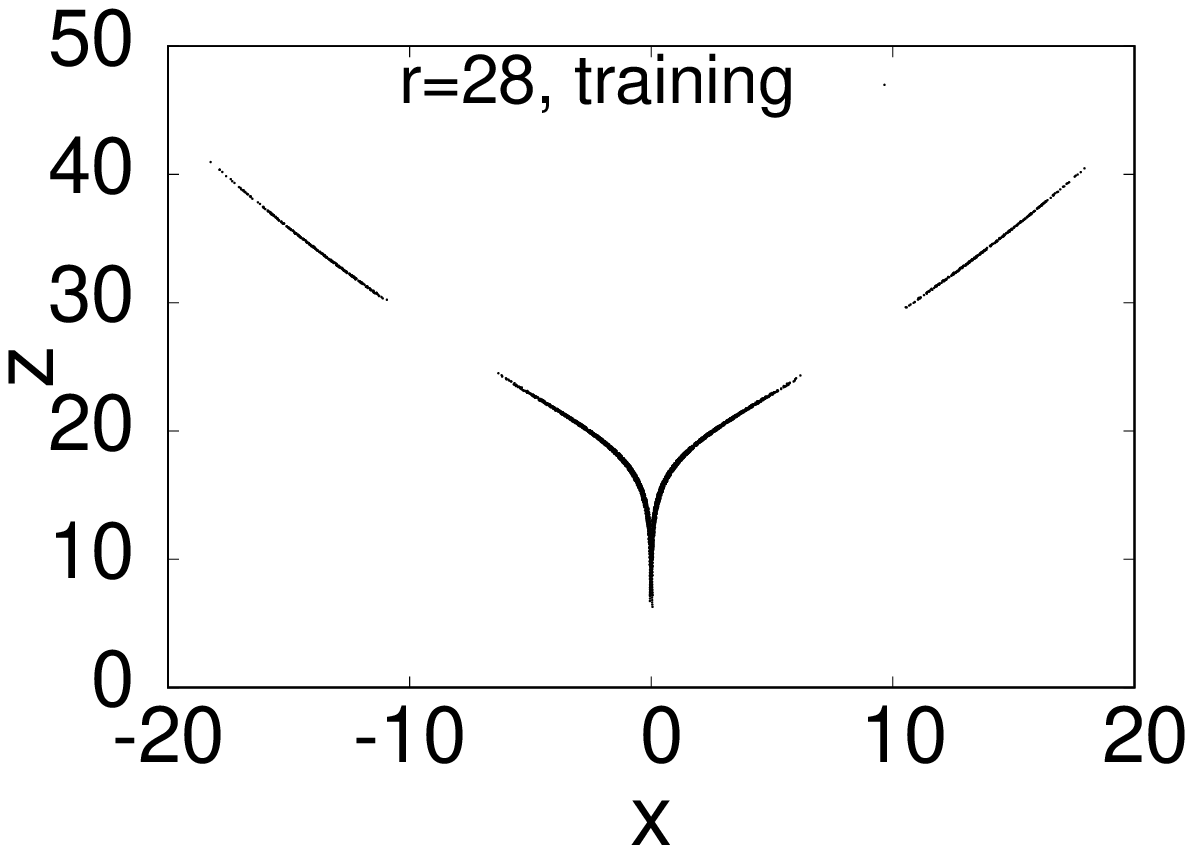}}%
		\subfigure[\hspace{-4mm}]{ \includegraphics[width=0.455\columnwidth,height=0.45\columnwidth]{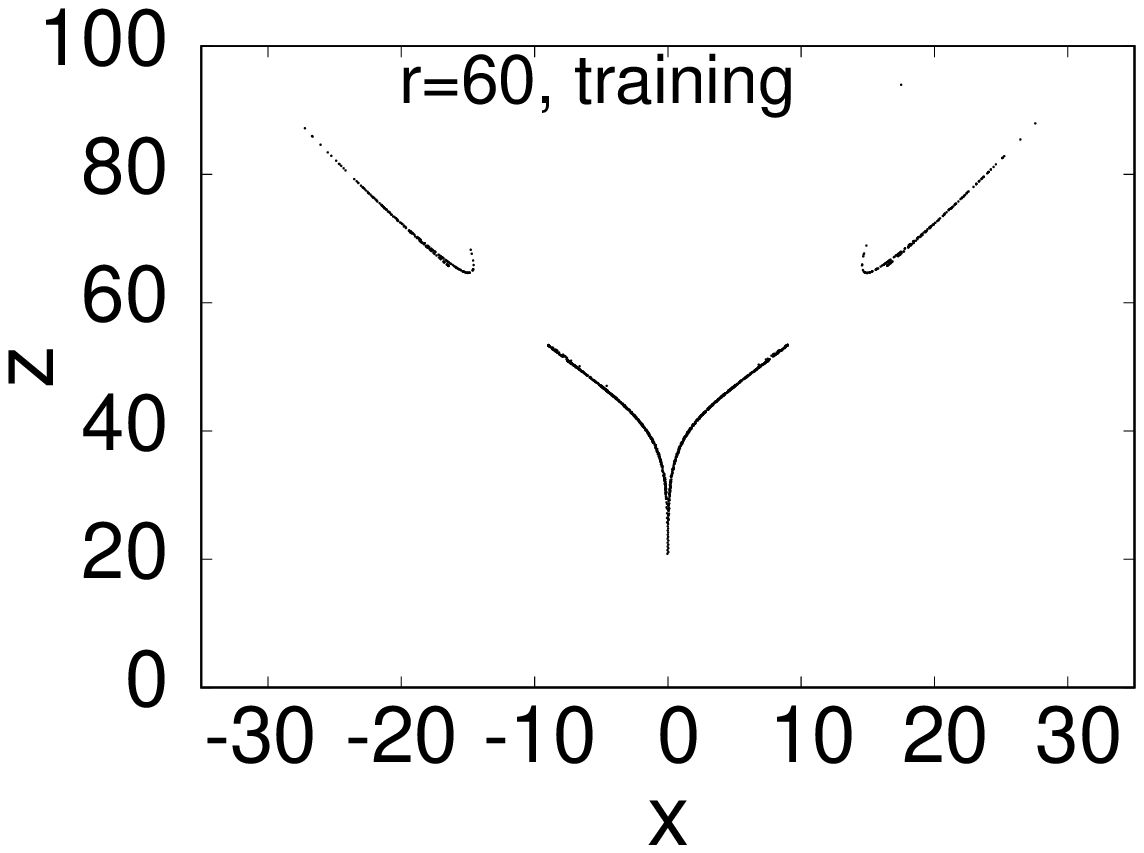}}%
\caption{{\bf Poincar\'e section-like plots}~($r=28$ (left) and $60$  (right)). The sets of points $(x,z)$: ((a), (b)) along a trajectory of the data-driven model using reservoir computing and ((c), (d)) along a long trajectory of the actual Lorenz system and ((e), (f)) along a short trajectory used for the training data are plotted when $|x-y|<\epsilon_p$, where $\epsilon_p=0.05$.
The time lengths of the three trajectories are $T=10^6$, $10^6$,
and $5000$, respectively. 
}
			\label{fig:poincare}
	\end{figure}

\begin{figure}
		\subfigure[\hspace{-2mm}]{\includegraphics[width=0.49\columnwidth,height=0.5\columnwidth]{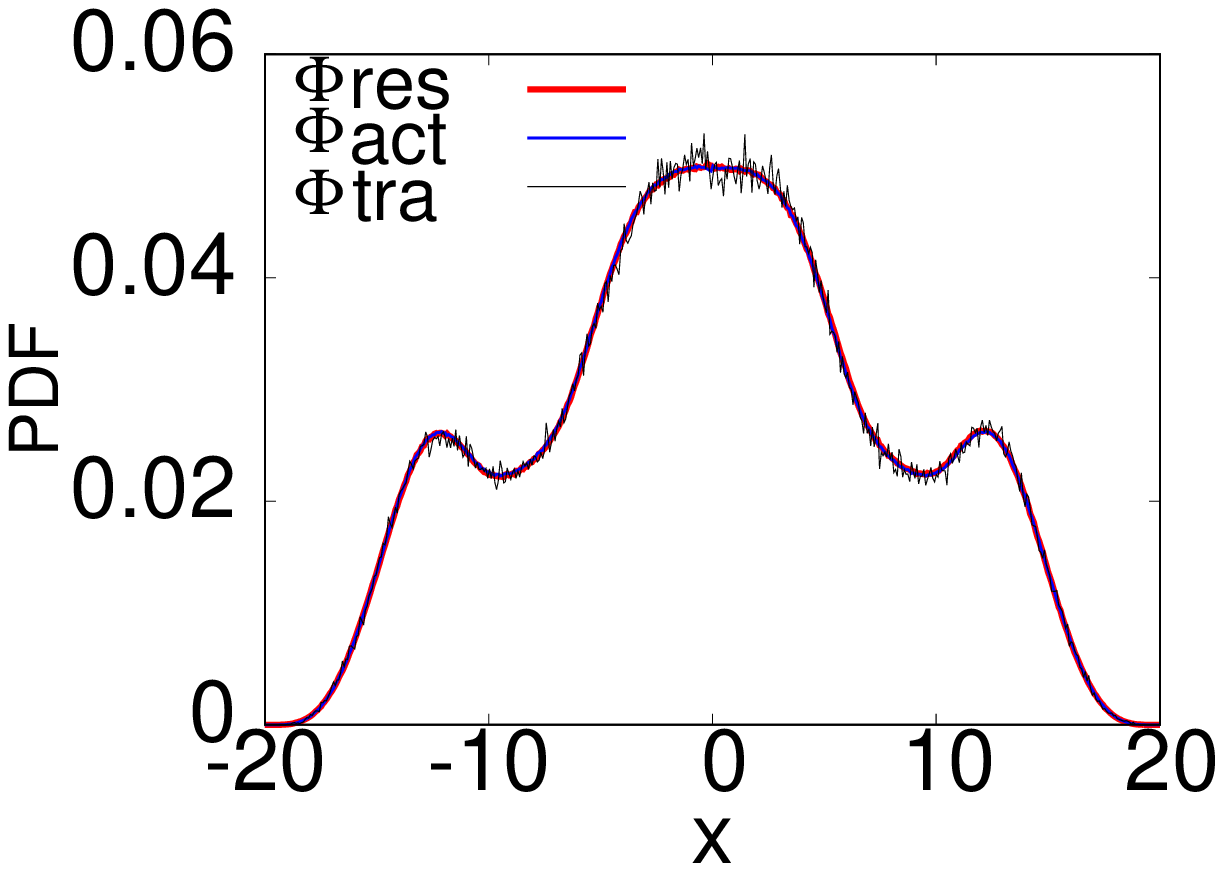}}%
		\subfigure[\hspace{-4mm}]{\includegraphics[width=0.49\columnwidth,height=0.5\columnwidth]{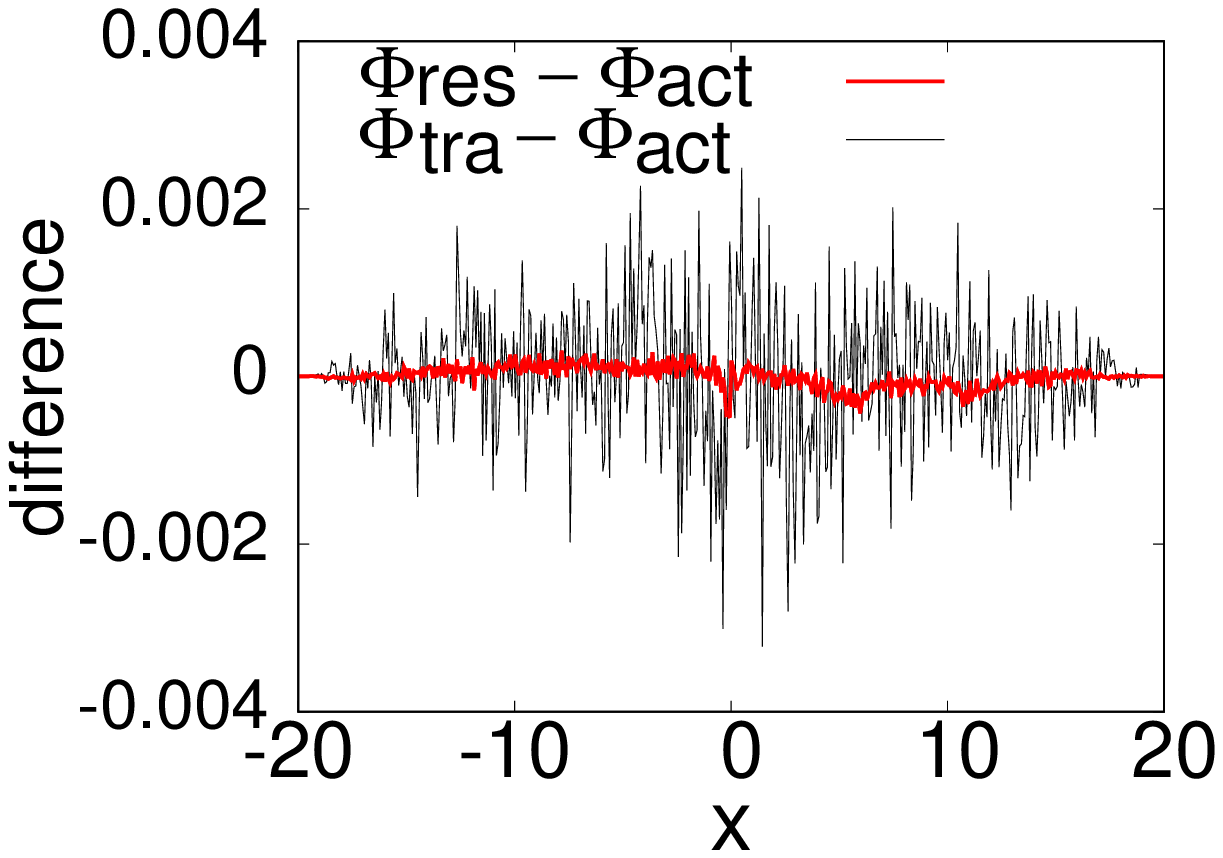}}\\%
		\subfigure[\hspace{-2mm}]{\includegraphics[width=0.49\columnwidth,height=0.5\columnwidth]{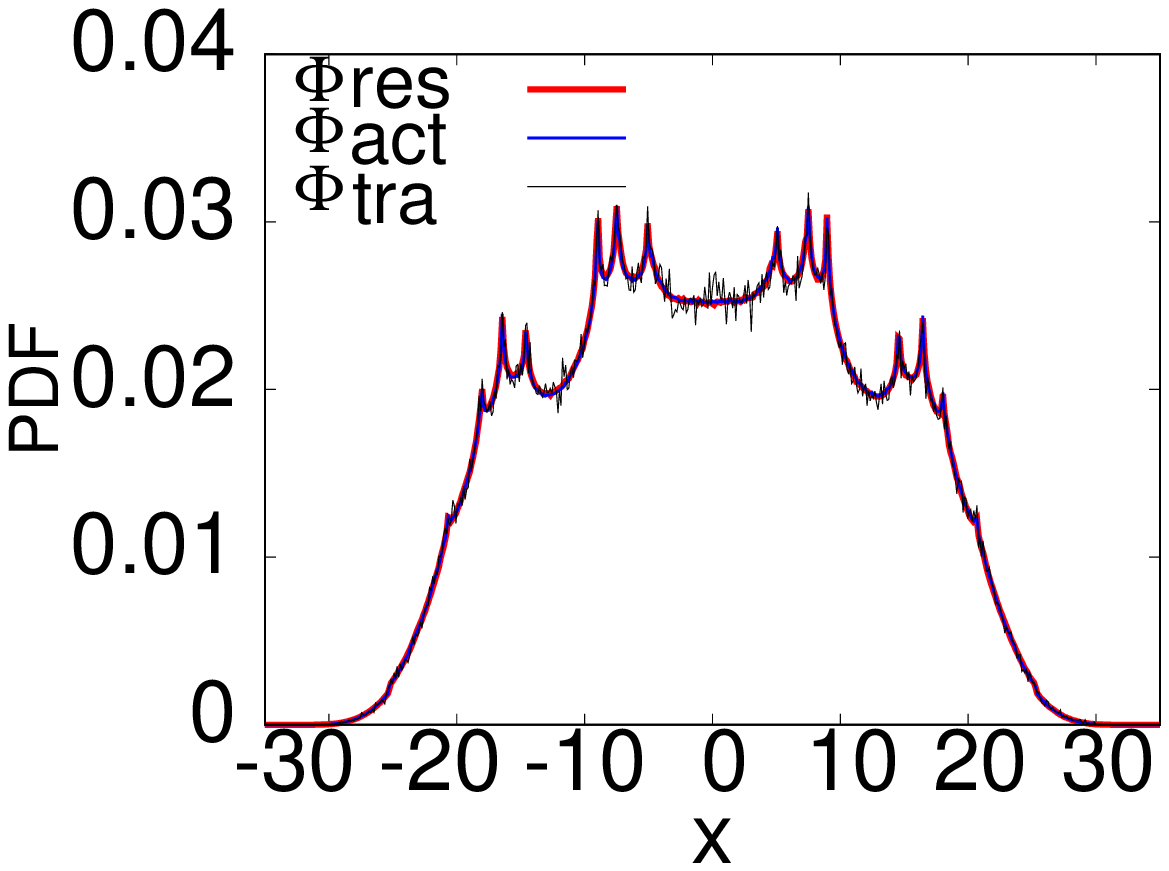}}%
		\subfigure[\hspace{-4mm}]{\includegraphics[width=0.49\columnwidth,height=0.5\columnwidth]{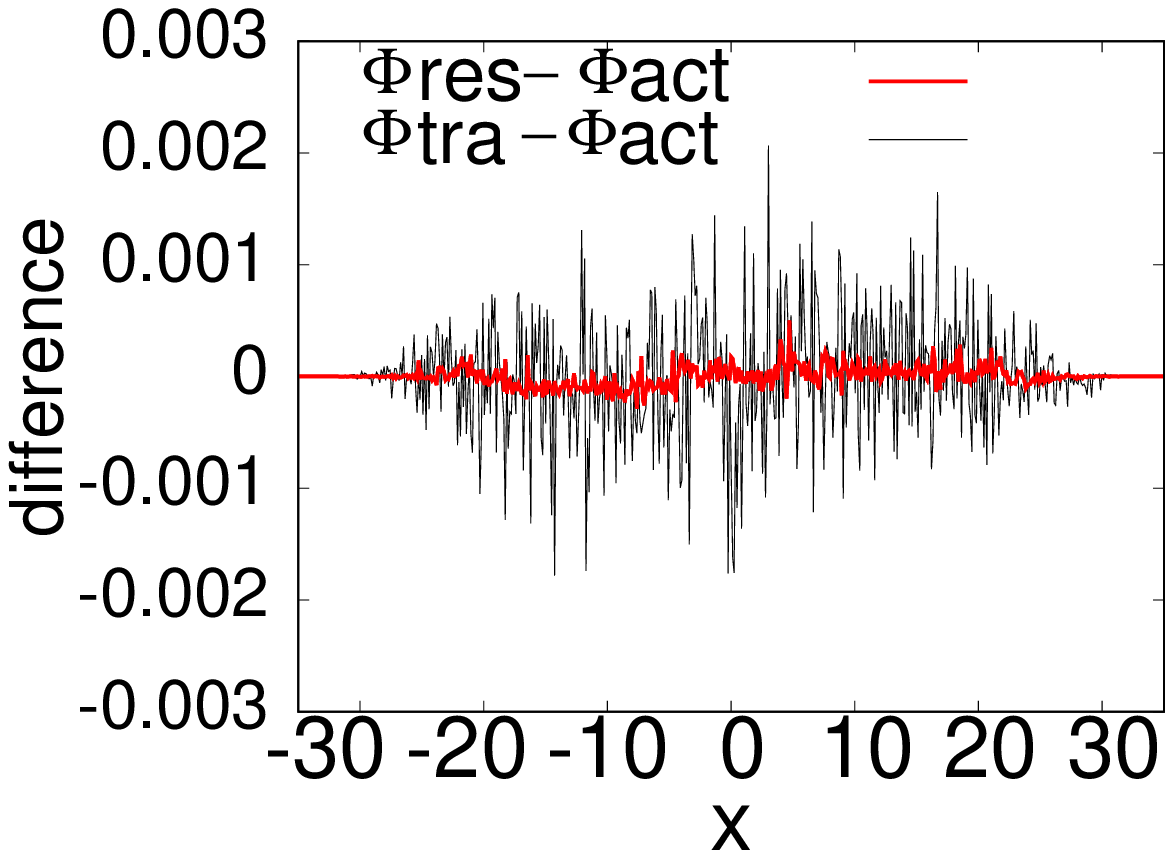}}
\caption{{\bf Density distributions of a variable}~($r=28$~(top) and $r=60$~(bottom)).
((a), (c)) The density distribution of the $x$ variable calculated from a length $T=10^6$ trajectory of the data-driven model by reservoir computing ($\Phi\text{res}$) is plotted together with that computed from the length $T=10^6$ long trajectory of the actual Lorenz system
($\Phi\text{act}$)
and 
with that computed from
the length $T=5000$ short trajectory 
($\Phi\text{tra}$) 
used as training data for constructing the data-driven model.
Here the length $T=10^6$ long trajectory is 
used to obtain a distribution which approximates 
the limiting distribution.
((b), (d)) The differences in the distributions are shown.
The average $\overline{x}$ and the standard deviation $\sigma$ of the density distribution 
 are as follows; $(\overline{x},\sigma)=(-0.004,7.924)$ for the data-driven model,
and $(0.009,7.925)$ for the Lorenz system with $r=28$;
$(\overline{x},\sigma)=(-0.018,12.092)$ for the data-driven model,
and $(-0.018,12.091)$ for the Lorenz system with $r=60$.
$\int|\Phi\text{res}-\Phi\text{act}|dx/
\int|\Phi\text{tra}-\Phi\text{act}|dx
\approx 1/8$ for $r=28$, and $\approx 1/6$ for $r=60$. 
}\label{fig:xdistribution-r28-60}
	\end{figure}

\begin{figure}
		\subfigure[\hspace{-2mm}]{\includegraphics[width=0.485\columnwidth,height=0.5\columnwidth]{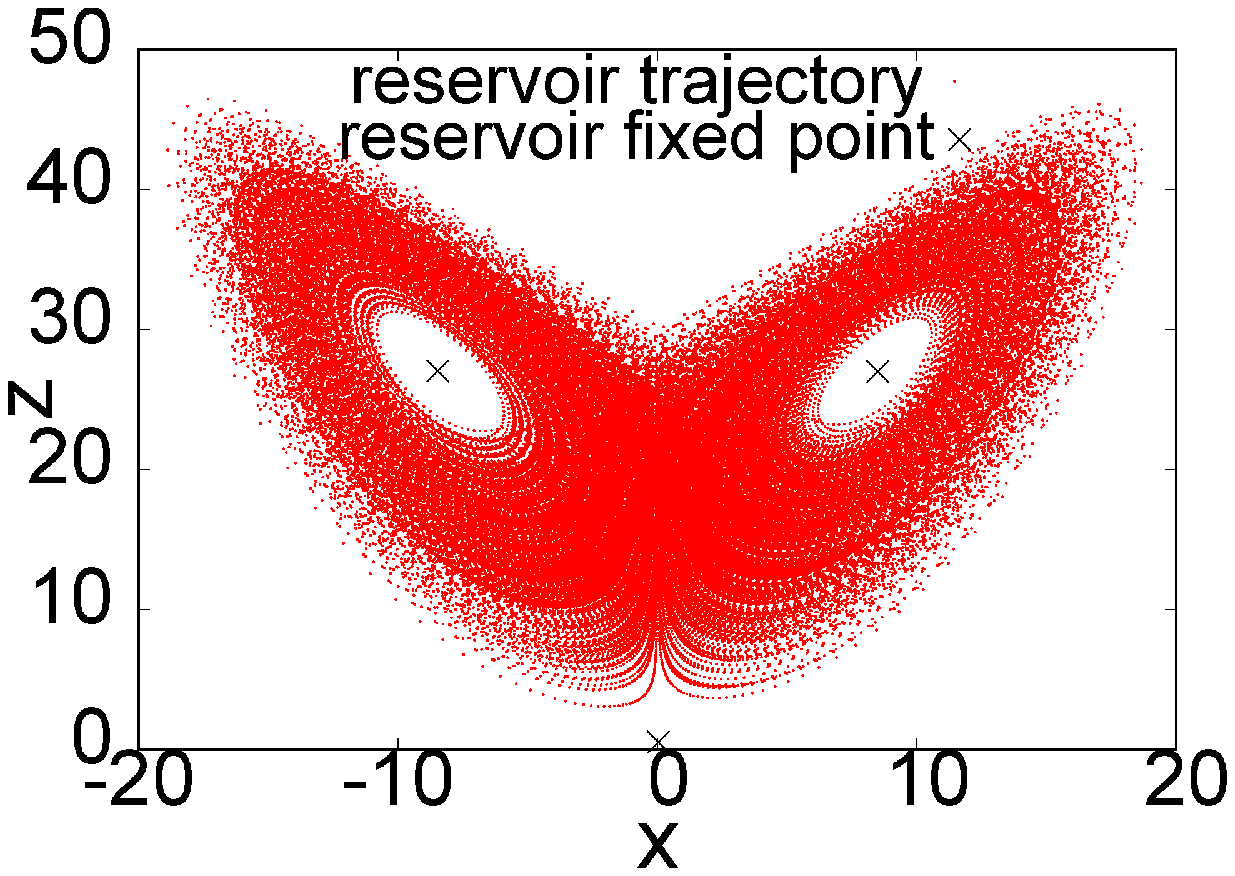}}%
		\subfigure[\hspace{-4mm}]{\includegraphics[width=0.485\columnwidth,height=0.5\columnwidth]{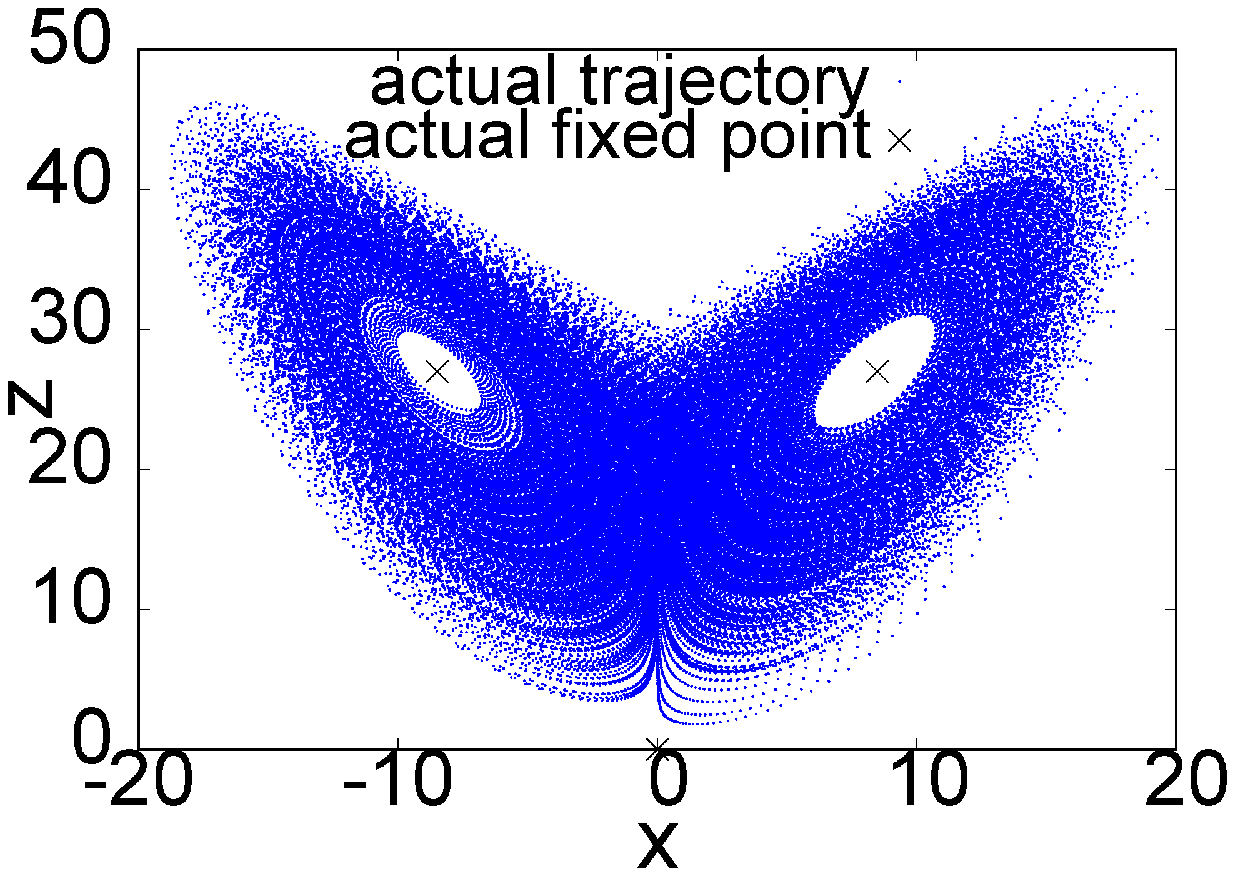}}
 		\caption{{\bf Fixed points.} 
(a) The three fixed points $(x^{*},y^{*},z^{*})$ of the data-driven model and (b) the corresponding unstable fixed points of the actual Lorenz system  are plotted together with trajectory points with 
the time length $T=10^4$.
The three fixed points of the  data-driven model are close to those of the actual Lorenz system, despite the fixed points being outside the training data, which is part of the actual trajectory.
See the coordinates and the eigenvalues of the Jacobian matrix at each fixed point in Table~\ref{tab:fixed-points}.
}\label{fig:fixedpoint}
	\end{figure}

\begin{figure}
		\subfigure[\hspace{-2mm}]{\includegraphics[width=0.45\columnwidth,height=0.55\columnwidth]{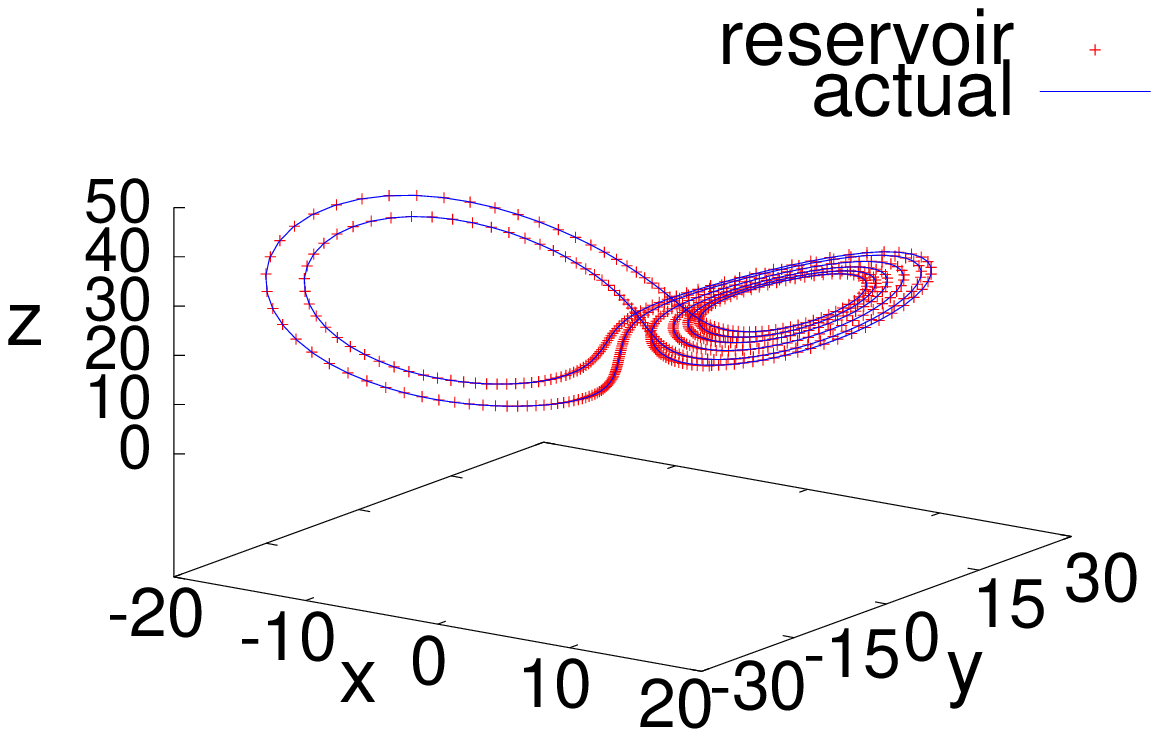}}%
		\subfigure[\hspace{-4mm}]{\includegraphics[width=0.535\columnwidth,height=0.55\columnwidth]{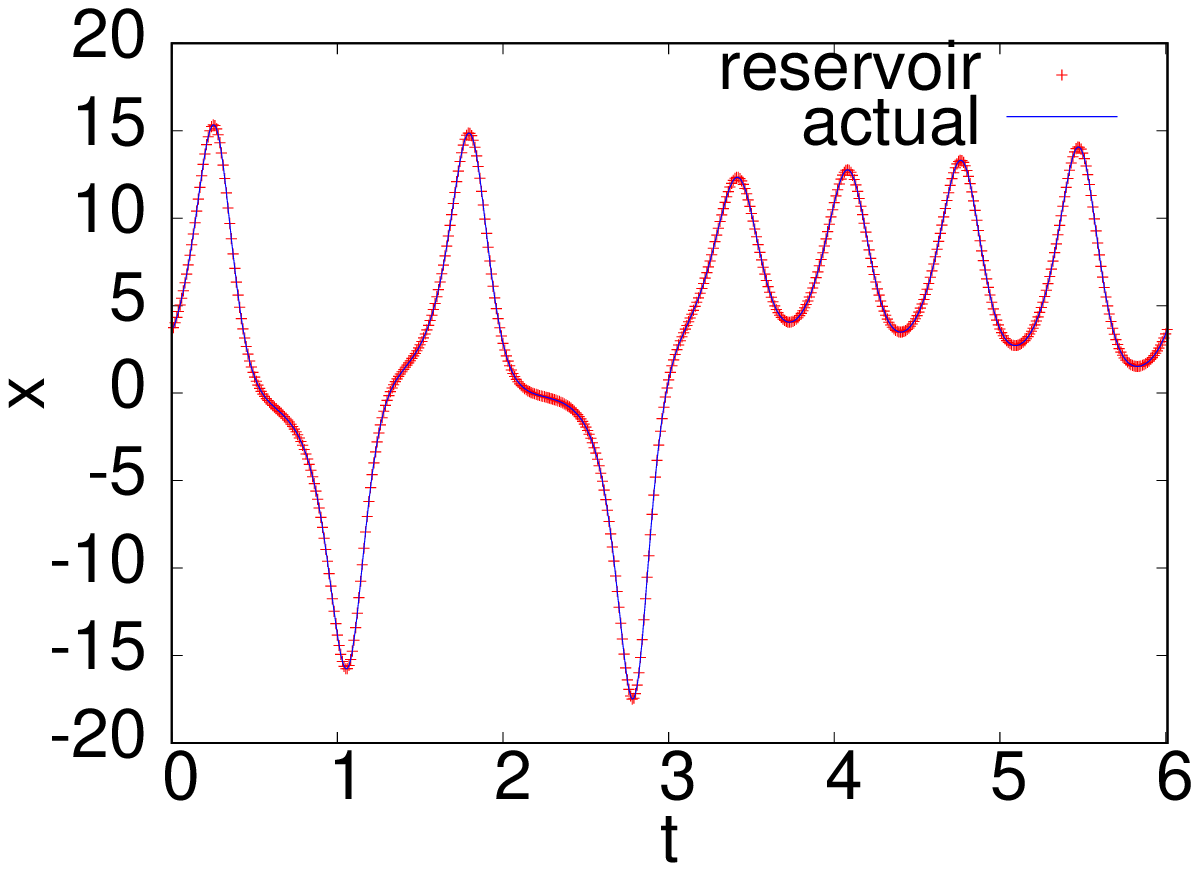}}
 		\caption{{\bf A periodic orbit-like trajectory.} 
 		(a) A periodic orbit-like trajectory obtained from the data-driven model is plotted 
 		together with the corresponding unstable  periodic orbit (period $T_p$ =$5.9973192969$) of the actual Lorenz system with $r=28$, and 
 		(b) their time developments of the $x$ variable. 
}
\label{fig:periodicorbit}
	\end{figure}

\indent{\bf Poincar\'e section-like plots.} 
The Poincar\'e section of the data-driven model of the Lorenz system has been studied~\cite{Pathak_2017}.
We compare the shape and size of the attractor of a data-driven model~(\ref{eq:reservoir}) with 
those of the attractor of the actual Lorenz system~(\ref{eq:lorenz}),  
and also with those of 
the set of points along 
the training time series data.
Figure~\ref{fig:poincare} presents their Poincar\'e section-like plots for $r=28$ and $60$.
For each of the two parameter cases, a set of trajectory points generated from the data-driven model seem to coincide with 
the chaotic attractor of the actual Lorenz system. 
Furthermore, the data-driven model has an attractor which is significantly larger than  the set of training data used to construct the model.\\
\indent{\bf Density distribution.} 
The density distribution of $x$ variable along a trajectory 
of the data-driven model is presented in Fig.~\ref{fig:xdistribution-r28-60}.
We compare the distribution with that obtained from the trajectory of the actual Lorenz system~(\ref{eq:lorenz}) and that calculated directly from the training data. 
The distribution of the actual Lorenz system can be captured 
by employing the data-driven model.
Remarkably, the distribution with a singular structure~\cite{zoldi_1998} in $r=60$ can be recovered using the data-driven model.\\
%
{\bf Fixed points and their stabilities.} 
Fixed points, which are fundamental structures of dynamical systems, are examined.
We identify fixed points in the space of the output variables directly, even though they were identified through the fixed points in the space of the reservoir state vector using the directional fibers method~\cite{krishnagopal_2019}. 
We also study the 
stability of each  unstable fixed point in the space of output variables.
For the data-driven model 
we consider a point $\mb{x^*}=(x^*,y^*,z^*)$ as a fixed point, when the following condition is satisfied:
    $\delta=\max_{n\in [0,n_0]} \| \mb{x}^* - \psi_{\mb{x}^*}(n\Delta t)  \|_{l^2} < \epsilon_0$
for some $\epsilon_0$ sufficiently small and for some $n_0$ sufficiently large, where $\psi_{\mb{x}^*}(n\Delta t)$ is the point iterated $n$ times from $\mb{x}^*$ by the data-driven model with the time step $\Delta t$. 
For the computation of a trajectory from a given point $\mb{x}^*$ of the data-driven model, reservoir state vector $\mb{r}(0)$ is determined to correspond to $\mb{x}^*$ by the pre-iterates. 
The echo state property~\cite{Jaeger_2001} in which our choice of parameters in the data-driven model~(\ref{eq:reservoir})
is satisfied guarantees that for 
each $\mb{x}^*$, the corresponding reservoir
state vector is determined uniquely.\\
\indent Table~\ref{tab:fixed-points} lists 
the obtained coordinates of the three fixed points,
$L_\text{res}$, $R_\text{res}$ and  $O_\text{res}$, together with those of the actual Lorenz system. 
We fix $(\epsilon_0,n_0)=(0.01,10000)$ 
for $L_{\text{res}}$ and $R_{\text{res}}$, and $(\epsilon_0,n_0)=(1,30)$ 
for $O_{\text{res}}$.
Figure~\ref{fig:fixedpoint} shows the fixed points 
together with the trajectory points.
Table~\ref{tab:fixed-points} also lists 
the eigenvalues of the Jacobian matrix at each fixed point. 
The values are obtained from the estimated formula of the Jacobian matrix described later for calculating the Lyapunov exponents and vectors.\\ 
\begin{table*}[t]
		\begin{tabular}{|l|r|r|r|r|r|r|r} 
			\hline  	
           	  & $L_{\text{res}}$& $R_{\text{res}}$& $O_{\text{res}}$& $L_{\text{actual}}$& $R_{\text{actual}}$& $O_{\text{actual}}$ \\ \hline
  			$x^*$ &$-8.47$ & $8.50$ & $0.04$ & $-8.49$ & $8.49$ & $0.00$  \\ \hline
  			$y^*$ &$-8.47$ & $8.50$ & $0.02$ & $-8.49$ & $8.49$ & $0.00$  \\ \hline
  			$z^*$ &$27.04$ & $27.01$ & $0.54$  & $27.00$ &$27.00$&$0.00$ \\ \hline
  			$\Lambda_1$&$0.09+10.19i$&$0.10+10.21i$&$11.67$&$0.09+10.20i$&$0.09+10.20i$&$11.83$\\ \hline
  			$\Lambda_2$&$0.09-10.19i$&$0.10-10.21i$&$-2.66$&$0.09-10.20i$&$0.09-10.20i$&$-2.67$\\ \hline
  			$\Lambda_3$&$-13.84$&$-13.86$&$-22.68$&$-13.85$&$-13.85$&$-22.83$\\ \hline
 \end{tabular}
 		\caption{{\bf Coordinates and eigenvalues of the Jacobian matrix at each of the three  unstable fixed points.} 
$L_{\text{res}}$, $R_{\text{res}}$,  and $O_{\text{res}}$ are fixed points of the data-driven model, whereas
$L_{\text{actual}}$, $R_{\text{actual}}$, and $O_{\text{actual}}$ are fixed points of the actual Lorenz system
with $r=28$. 
The coordinates $(x^{*},y^{*},z^{*})$ and the
eigenvalues $(\Lambda_1,\Lambda_2,\Lambda_3)$ of the Jacobian matrix at each fixed point of the data-driven model are close to those of the corresponding fixed point of the actual Lorenz system.
}
 		
 		\label{tab:fixed-points}
\end{table*}
%
\indent{\bf Periodic trajectory.} 
Periodic orbits are also the fundamental structures of dynamical systems. 
We confirm that the data-driven model of discrete time has a periodic orbit-like trajectory that travels near the corresponding periodic orbit of the actual Lorenz system~(\ref{eq:lorenz}) of continuous time.
We call $\{\psi_{\mb{x}(0)}(n\Delta t)\}_{n\in [0,n_p]}$  a periodic orbit-like trajectory, 
if the following value
is sufficiently small for a periodic trajectory $\{\mb{x}(t)\}$ of period $T_p$ of the actual Lorenz system:
    $\delta_p=\max_{n\in [0,n_p]} \| \mb{x}(n\Delta t) - \psi_{\mb{x}(0)}(n\Delta t)\|_{l^2},$
where 
$n_p$ is the smallest integer satisfying 
$n_p\Delta t\ge T_p$.
Among the periodic orbit-like trajectories of the data-driven model corresponding to the 50 periodic orbits with low periods, $\delta_p <0.1$ for 40 cases and $\delta_p <0.4$ for the other 10 cases. 
Figure~\ref{fig:periodicorbit} gives an example of a periodic orbit-like trajectory, which has the largest value of $\delta_p$ among the 50 periodic orbits with low periods.\\
%
%
\indent{\bf Lyapunov exponents and Lyapunov vectors.}  
\begin{table*}
 	\begin{tabular}{|l|r|r|r|r|r|r|r|r|}
			\hline  	
           	  $r$ & $\lambda^{(1)}_{\text{res}}$& $\lambda^{(2)}_{\text{res}}$& $\lambda^{(3)}_{\text{res}}$& $D^{KY}_{\text{res}}$& $\lambda^{(1)}_{\text{actual}}$& $\lambda^{(2)}_{\text{actual}}$& $\lambda^{(3)}_{\text{actual}}$& $D^{KY}_{\text{actual}}$\\ \hline
  			$28$ & $0.901$ & $0.000$ & $-14.570$ & $2.06$ & $0.902$ & $0.000$ & $-14.570$ & $2.06$ \\ \hline
  			$60$ & $1.402$ & $0.000$ & $-15.070$ & $2.09$ & $1.404$ & $0.000$ & $-15.071$ & $2.09$ \\ \hline
 \end{tabular}
 		\caption{{\bf Lyapunov exponents and Lyapunov dimensions.} Lyapunov exponents  of the data-driven model using reservoir computing $(\lambda^{(1)}_\text{res},\lambda^{(2)}_\text{res},\lambda^{(3)}_\text{res})$ 
 		 and those of the actual Lorenz system $(\lambda^{(1)}_\text{actual},\lambda^{(2)}_\text{actual},\lambda^{(3)}_\text{actual})$ are listed.
 		The values are computed using the 
 		 four-stage and fourth-order Runge--Kutta
 		method with time step $2\Delta t$ from the points along an orbit trajectory and the estimated Jacobian matrices.
 		 The Lyapunov dimensions 
 		 $D^{KY}_\text{res}$ for the data-driven model and $D^{KY}_\text{actual}$ for the actual Lorenz system 
 		  are estimated from the Kaplan--Yorke formula~\cite{kaplan_1979}.
}
 		\label{tab:lyapunov-exponents}
\end{table*}
The Lyapunov exponents are used to evaluate the degree of instability and  
estimate the Lyapunov dimension of a dynamical system.
In some studies, the Lyapunov exponents 
of a data-driven model by reservoir computing
were calculated in the space of $N$-dimensional reservoir state vector~\cite{Pathak_2017,Pathak_2018, gallicchio_2018,pyragas_2020}. 
\textcite{Pathak_2017} computed Lyapunov exponents for the 
reservoir state vector and found that they almost coincide with 
those of the original system for the case of a partial differential equation, whereas only positive and neutral exponents coincide with those for the Lorenz system.
To the best of the authors' knowledge, 
they have not been computed in a space of output variables.

First we compute the first Lyapunov exponent using the traditional method which has been used to estimate the Lyapunov exponent from an experimental data without the knowledge of the equation~\cite{wolf_1985}. 
The first Lyapunov exponent estimated from a time-series of the data-driven model and that of the actual Lorenz system 
are 0.962 and 0.954, respectively\footnote{We choose  parameters in \cite{wolf_1985} to be $(\text{DIM}, \text{TAU}, \text{SCALMX}, \text{SCALMN}, \text{EVOLV}, \text{ANGLMX}_{\text{main}})$$=(3,11,0.1,0.001,600,0.013)$. Note that the estimated exponents are found to be robust (within $10\%$ of the error) under the choices of parameters $(\text{SCALMX}, \text{EVOLV}, \text{ANGLMX}_{\text{main}})=(0.1\pm0.01,600\pm20,0.013\pm0.002)$.}.

\begin{figure}
		\subfigure[\hspace{-2mm}]{\includegraphics[width=0.485\columnwidth,height=0.5\columnwidth]{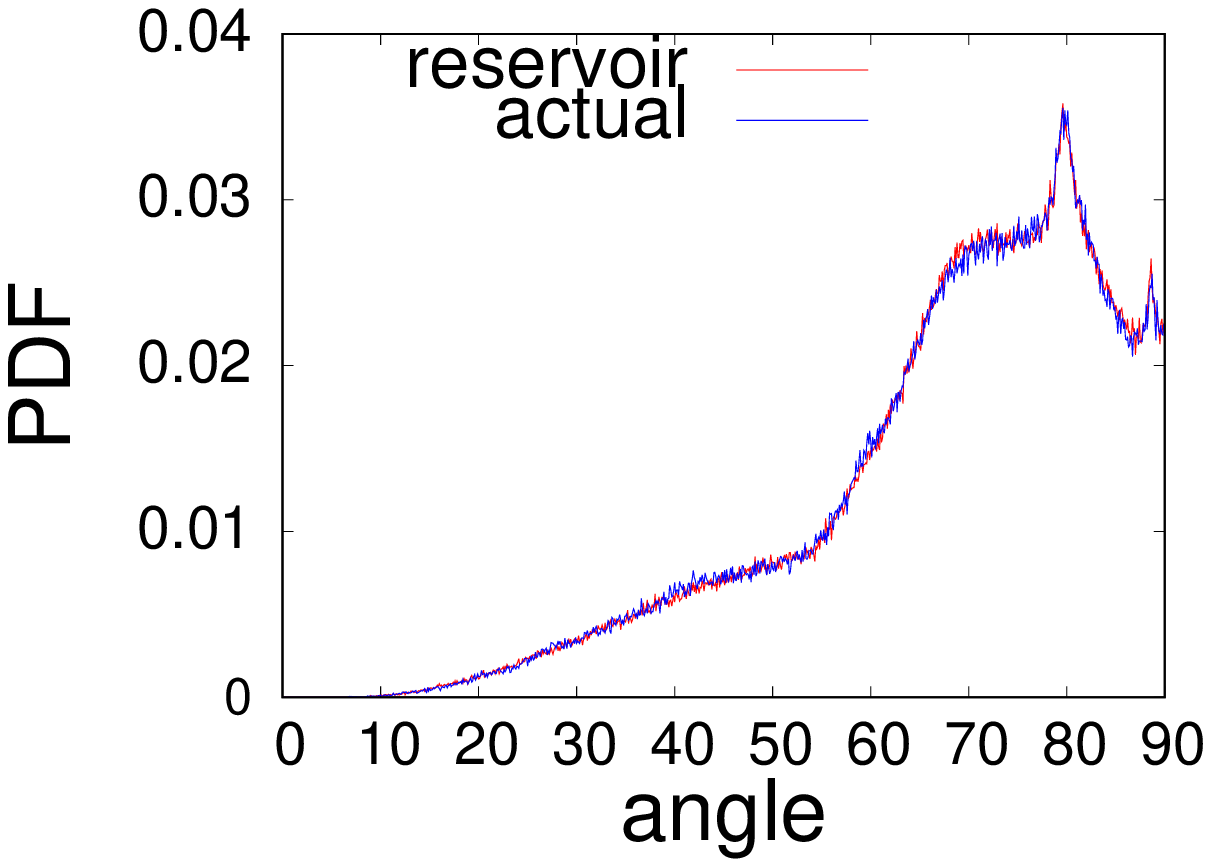}}%
		\subfigure[\hspace{-4mm}]{\includegraphics[width=0.485\columnwidth,height=0.5\columnwidth]{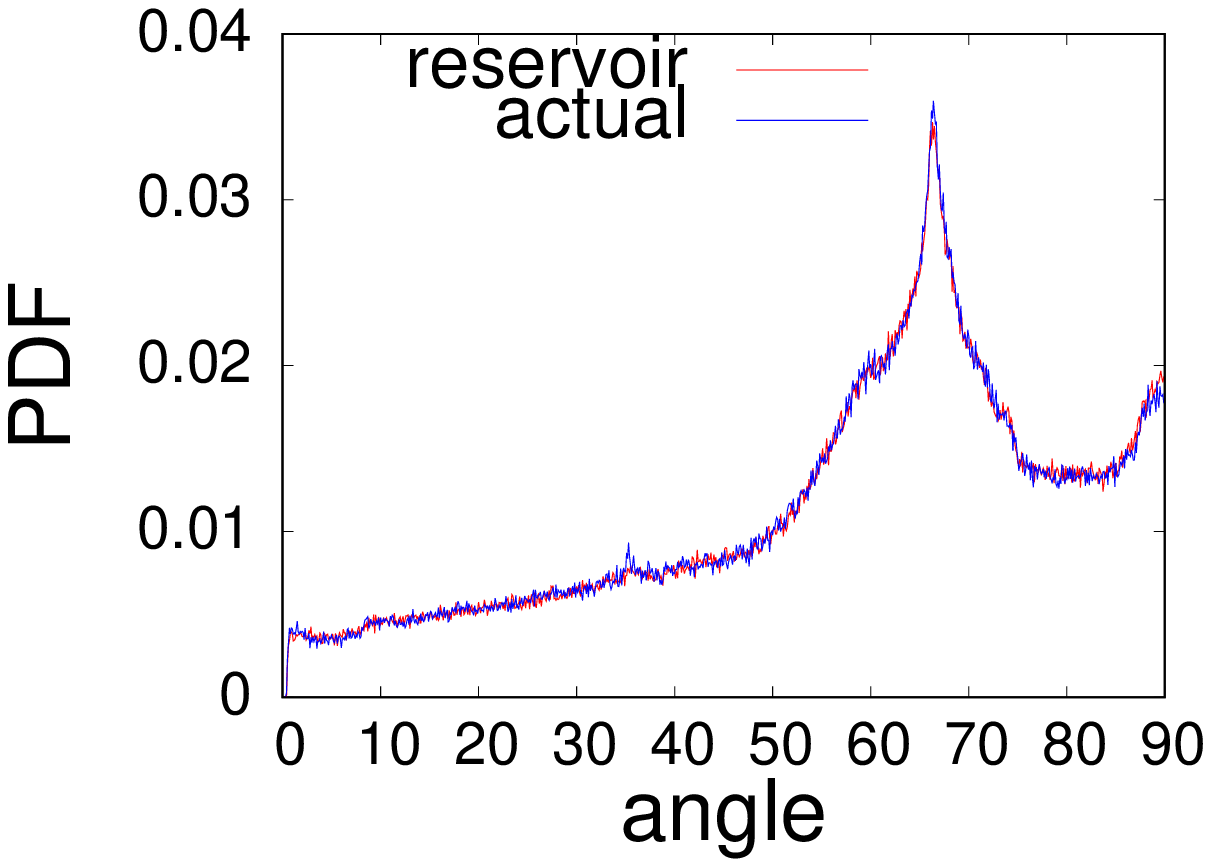}}
    \caption{\label{fig:angle-dist}
    {\bf Distribution of the angle between stable and unstable manifolds 
    along a trajectory 
    }~((a) $r=28$ and (b) $r=60$).
      The density distribution of the manifold angles (degree) at points along a trajectory is shown for a data-driven model using reservoir computing  together with that of the actual Lorenz system.
      }
	\end{figure}

In this paper, an attempt is made to compute Lyapunov exponents in the space of output variables corresponding to $x,y$ and $z$ for the Lorenz system.
 Here we describe how to compute Lyapunov exponents and vectors in the original variables numerically from a trajectory of the data-driven model. 
We first estimate the Jacobian matrix at each point $(x,y,z)$ along the trajectory of the data-driven model as follows: 
(i) Apply the Taylor series expansion of order six to estimate $\dot{x}=dx/dt$, $\dot{y}=dy/dt$ and $\dot{z}=dz/dt$ at each sample point along the discrete trajectory; 
(ii) Apply linear regression to the estimated values of $\dot{x}$, $\dot{y}$ and $\dot{z}$ 
by $x^l y^m z^n$ $(0\le l+m+n \le 3,~l,m,n\ge 0)$ as explanatory  variables; 
(iii) Obtain the Jacobian matrix $J({\bf x})$ at each point ${\bf x}$ by differentiating polynomials with the regression coefficients estimated in (ii).\\
\indent We compute Lyapunov exponents and vectors by integrating the 
linear ordinary differential equation having coefficients 
determined by the Jacobian matrices ($\dot{\bf x}(t)=J({\bf x}(t)){\bf x}(t)$), while the orbit 
is given by the trajectory of the data-driven model.
Note that in this computation the discrete time trajectory points of a data-driven model are considered samples of the continuous time trajectory. 
For the high-accuracy computation with a rather large 
time step $\Delta t=0.05$ of the reservoir computing, 
we employ four-stage and fourth-order Runge--Kutta method with 
time step $2\Delta t$ from the points along an orbit trajectory.\\
\indent The results are compared with those of the actual Lorenz system~(\ref{eq:lorenz}) for two sets of parameters. 
Table~\ref{tab:lyapunov-exponents} shows the
agreement of the Lyapunov exponents and the Lyapunov dimensions.\\
We also compute (co-variant) Lyapunov vectors,
which measure the degree of hyperbolicity by  calculating the angle between the stable and unstable manifolds at some trajectory point~\cite{ginelli_2007}. 
  \\
  \indent{\bf Manifold structure and Tangency.}
Using the computed Lyapunov vectors we investigate the manifold structures of the data-driven model, particularly the degree of hyperbolicity and the tangencies between the stable and the unstable manifolds. 
We consider the Lorenz system of $r=28$ without tangencies and of $r=60$ with tangencies for the comparison~\cite{saiki_2010}. 
Figure~{\ref{fig:angle-dist}} shows the probability density function 
of an angle between a tangent vector of a 
stable manifold and that of an unstable manifold along an orbit trajectory for each of the actual system and the data-driven model.
For each case of the parameters, $r=28$ and $r=60$, the angle distributions are quite similar in shape, indicating that the data-driven model can reconstruct the manifold structures.
Moreover, Fig.~{\ref{fig:angle-dist}}~(b) suggests that the data-driven model can represent a non-hyperbolic structure with tangencies between stable and unstable manifolds.\\ 

%
%
%
%
%
%
%
%
%
%
\section{IV. R\"ossler system}
\begin{figure}
		\subfigure[\hspace{-2mm}]{\includegraphics[width=0.485\columnwidth,height=0.5\columnwidth]{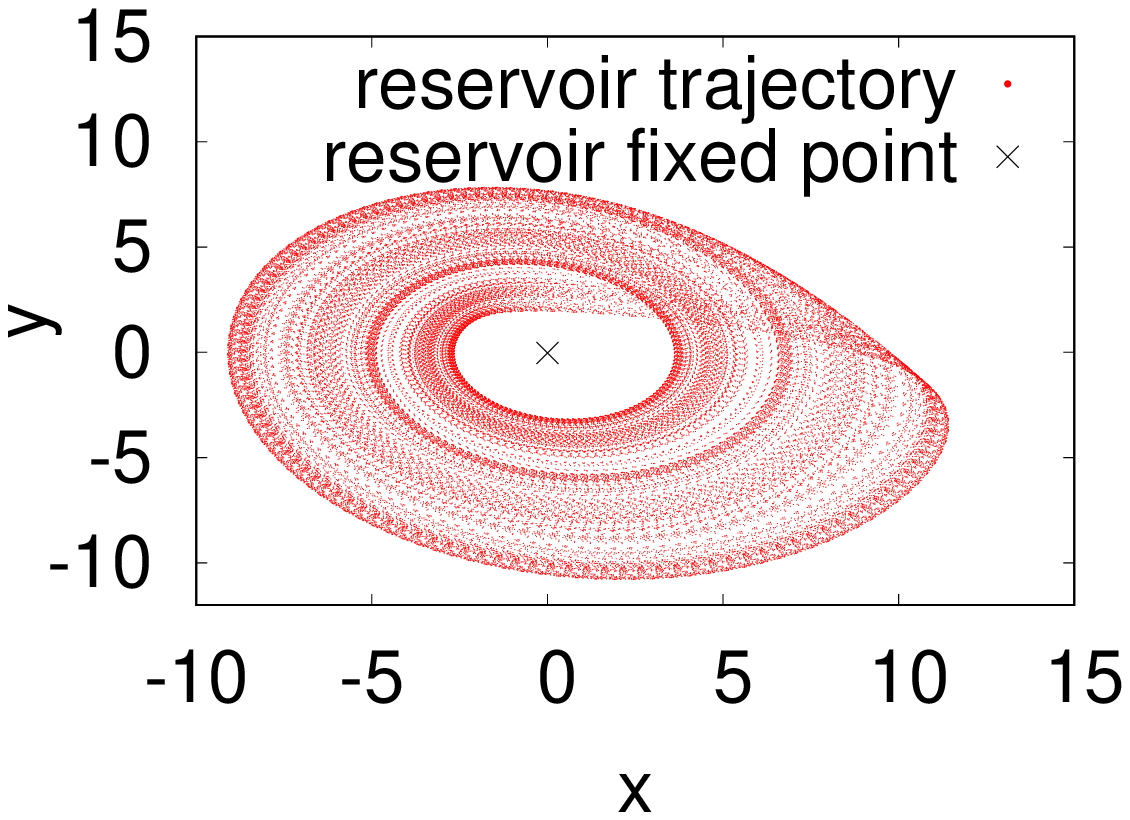}}%
		\subfigure[\hspace{-4mm}]{\includegraphics[width=0.485\columnwidth,height=0.5\columnwidth]{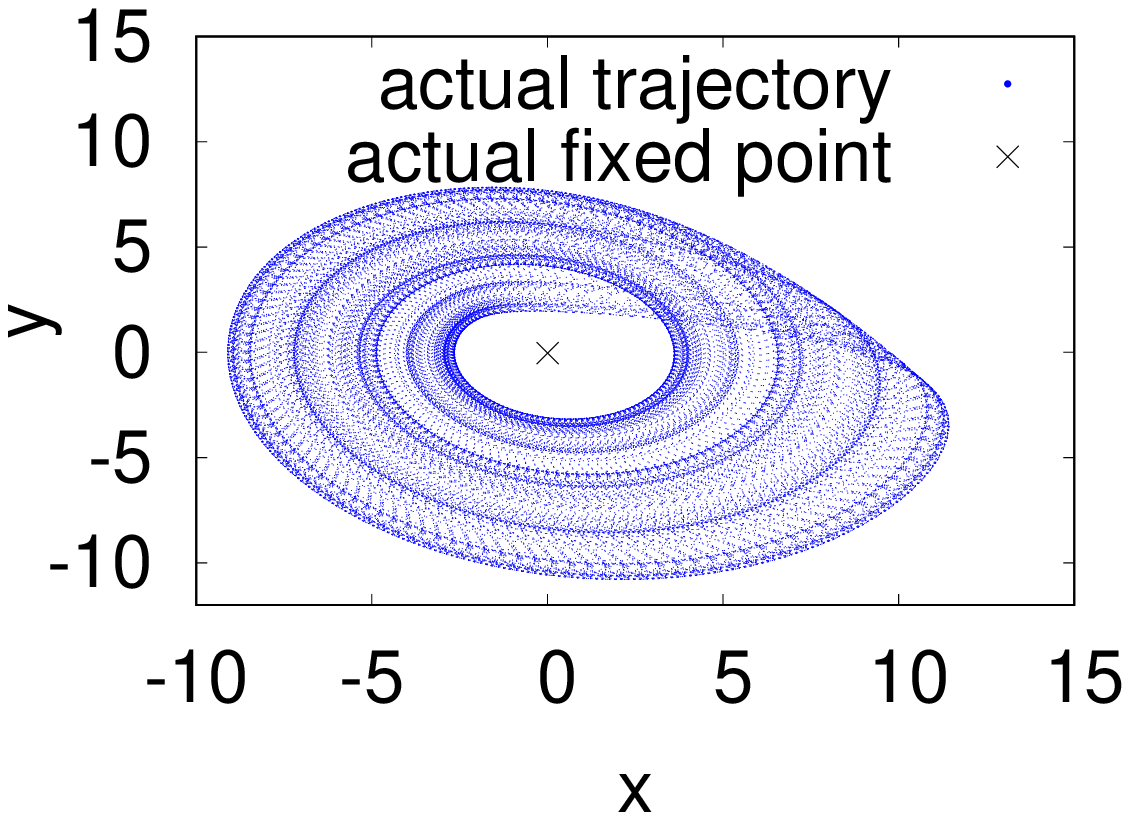}}
			\caption{{\bf 
			Fixed point.} 
			(a) The fixed point $(x^{*},y^{*},z^{*})$ of the data-driven model  and (b) the corresponding unstable fixed point of the actual R\"{o}ssler system are plotted together with trajectory points with the time length $T=2500$. 
			The fixed point of the data-driven model are close to those of the actual R\"{o}ssler system.  
			See the coordinates and the eigenvalues of the Jacobian matrix at each fixed point in Table~\ref{tab:fixed-points-rossler}.
			}
			\label{fig:fixedpoint-rossler}
	\end{figure}

\begin{figure}
		\subfigure[\hspace{-2mm}]{\includegraphics[width=0.485\columnwidth,height=0.5\columnwidth]{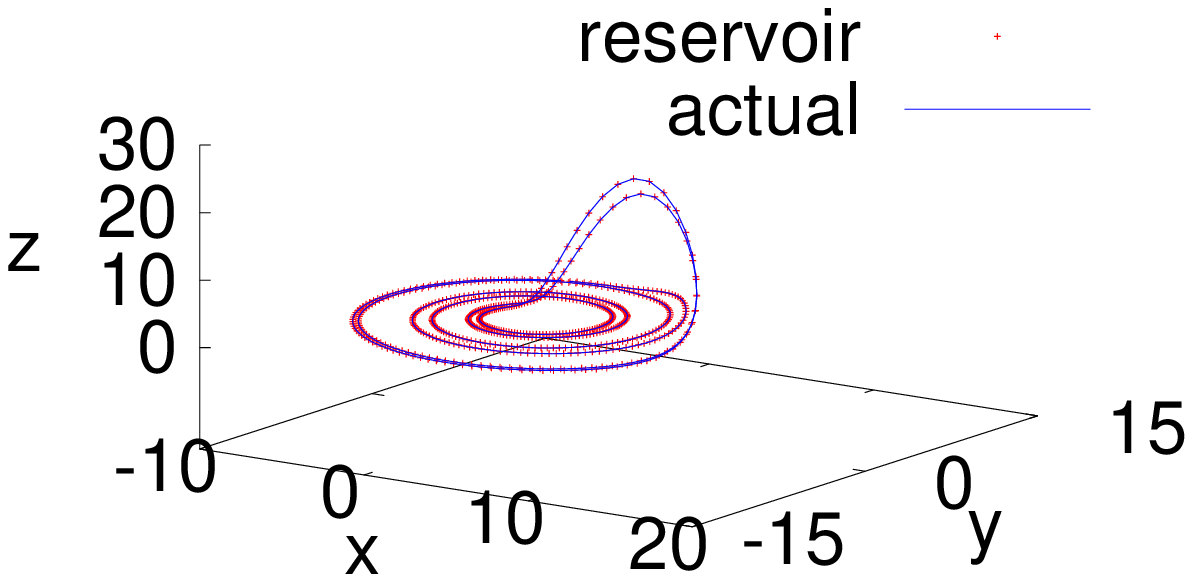}}%
		\subfigure[\hspace{-4mm}]{\includegraphics[width=0.485\columnwidth,height=0.5\columnwidth]{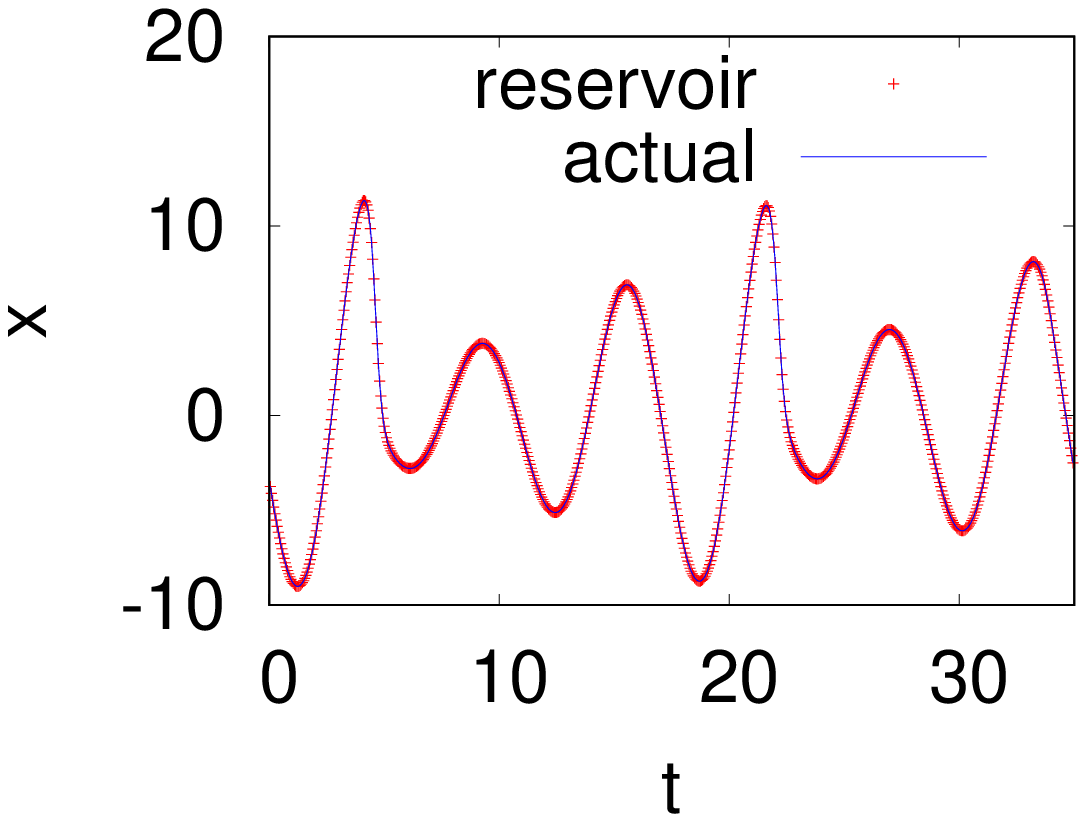}}
			\caption{{\bf 
			A periodic orbit-like trajectory.} 
     	A periodic orbit-like trajectory obtained from the data-driven model is plotted 
 		together with the corresponding unstable  periodic orbit with period
 		 $T_{p}$ =$35.06122601174815$.
 		 (a) The projections and (b) the time-series are shown.
		}
			\label{fig:periodicorbit-rossler}
	\end{figure}

\indent 
We confirm that for the R\"ossler system a data-driven model using reservoir computing has quite similar dynamical system properties to those of the original system.
\\
\indent
{\bf Fixed points and their stabilities.} 
\indent Table~\ref{tab:fixed-points-rossler} lists the obtained coordinates of a fixed point, $F_\text{res}$, together with that of the actual R\"ossler system. We fix $(\epsilon_0,n_0)=(0.01,800)$ for $F_{\text{res}}$. Figure~\ref{fig:fixedpoint-rossler} shows the fixed points together with the trajectory points. Table~\ref{tab:fixed-points-rossler} also lists the eigenvalues of the Jacobian matrix at the fixed point. The values are obtained from the estimated formula of the Jacobian matrix described later for calculating the Lyapunov exponents and vectors. 
%
\begin{table*}
		\begin{tabular}{|l|r|r|r|r|r|r|r} 
			\hline  	
           	  	&  $x^*$ &$y^*$ &$z^*$ &$\Lambda_1$&$\Lambda_2$&$\Lambda_3$  \\ \hline
$F_{\text{actual}}$&$0.0070$&$-0.0351$&$0.0351$&$0.0970 + 0.9952 i$&$0.0970 - 0.9952 i$&$-5.6870$  \\ \hline
$F_{\text{res}}$&$0.0015$&$-0.0315$&$ 0.0317$&$0.0926 + 0.9702 i$&$0.0926 - 0.9702 i$&$-5.6833$\\ \hline
 \end{tabular}
 		\caption{{\bf 
 		Coordinates and eigenvalues of the Jacobian matrix at each of the three  unstable fixed points.} 
$F_{\text{res}}$ are fixed points of the data-driven model, whereas
$F_{\text{actual}}$ are fixed points of the actual R\"{o}ssler system. 
The coordinates $(x^{*},y^{*},z^{*})$ and the
eigenvalues $(\Lambda_1,\Lambda_2,\Lambda_3)$ of the Jacobian matrix at each fixed point of the data-driven model are close to those of the corresponding fixed point of the actual R\"{o}ssler system.
}
			\label{tab:fixed-points-rossler}
\end{table*}
\\
\indent
{\bf Periodic trajectory.} 
We confirm that the data-driven model of discrete time has a periodic orbit-like trajectory that travels near the corresponding periodic orbit of the actual 
R\"{o}ssler system~(\ref{eq:rossler}) of continuous time. 
Figure~\ref{fig:periodicorbit-rossler} gives an example of 
{periodic orbit-like trajectories}.\\
%

%
{\bf Lyapunov exponents and Lyapunov vectors.}  
We compute the Lyapunov exponents in the space of output variables corresponding to $x,y$ and $z$ for the R\"ossler system. 
The results are compared with those of the actual R\"ossler system~(\ref{eq:rossler}). 

By using the Wolf's method we compute the first Lyapunov exponent from a  time-series of the data-driven model as $0.0707$ and that of the actual system as $0.0708$, which almost coincide with each other.

\begin{figure}
    \begin{center}
      \begin{tabular}{c}
	     \includegraphics[width=0.485\columnwidth,height=0.5\columnwidth]{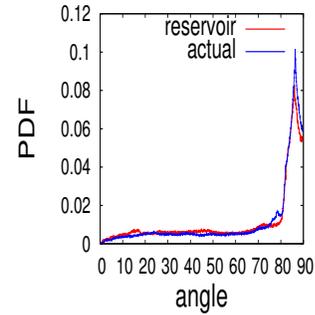}
      \end{tabular}
    \end{center}
    \caption{\label{fig:angle-dist-rossler}
    {\bf 
    Distribution of the angle between stable and unstable manifolds 
    along a trajectory of the R\"ossler system.
    }
      The density distribution of the manifold angles (degree) at points along a trajectory is shown for a data-driven model using reservoir computing  together with that of the actual R\"ossler system.
      }
  \end{figure}

Table~\ref{tab:lyapunov-exponents-rossler} shows the
agreement of the Lyapunov exponents and the Lyapunov dimensions by using our method.
We also compute (co-variant) Lyapunov vectors, which measure the degree of hyperbolicity by calculating the angle between the stable and unstable manifolds at some trajectory point.  
\begin{table*}
\vspace{10mm}
 	\begin{tabular}{|r|r|r|r|r|r|r|r|r|} 
			\hline  	
           	  $\lambda^{(1)}_{\text{res}}$& $\lambda^{(2)}_{\text{res}}$& $\lambda^{(3)}_{\text{res}}$& $D^{KY}_{\text{res}}$&  $\lambda^{(1)}_{\text{actual}}$& $\lambda^{(2)}_{\text{actual}}$& $\lambda^{(3)}_{\text{actual}}$& $D^{KY}_{\text{actual}}$\\ \hline
  			$0.07150$ & $0.00004$ & $-5.38813$ & $2.013$ & $0.07151$ & $0.00001$ & $-5.38809$  & $2.013$\\ \hline
 \end{tabular}
 		\caption{{\bf 
 		Lyapunov exponents and Lyapunov dimensions.} Lyapunov exponents  of the data-driven model using reservoir computing $(\lambda^{(1)}_\text{res},\lambda^{(2)}_\text{res},\lambda^{(3)}_\text{res})$ 
 		 and those of the actual R\"ossler system $(\lambda^{(1)}_\text{actual},\lambda^{(2)}_\text{actual},\lambda^{(3)}_\text{actual})$ are listed.
 		The values are computed using the 
 		 four-stage and fourth-order Runge--Kutta
 		method with time step $2\Delta t$ from the points along an orbit trajectory and the estimated Jacobian matrices.
}
			\label{tab:lyapunov-exponents-rossler}
\end{table*}
\\
\indent
{\bf Manifold structure and Tangency.}
The degree of hyperbolicity and the tangencies between the stable and the unstable manifolds are investigated for 
the R\"ossler system.
Figure~{\ref{fig:angle-dist-rossler}} shows the probability density function of an angle between a tangent vector of a stable manifold and that of an unstable manifold along an orbit trajectory for each of the actual system and the data-driven model. 
The angle distributions are quite similar in shape, indicating that the data-driven model can reconstruct the manifold structures. 

\section{V. chaotic fluid flow.}
\indent{\bf Laminar lasting time distribution of chaotic fluid flow.} 
We have clarified that hyperbolic fixed points and their eigenvalues are estimated in high accuracy by constructing a data-driven model using 
reservoir computing, even if the training time-series data are far away from the 
fixed points. 
We study a chaotic flow in macroscopic variables whose behavior has a random switching between laminar and bursting states.
Here we consider the set of laminar state as a certain chaotic saddle and compute the lasting time distribution staying in the neighborhood of it. 
We are interested in the lasting time distribution 
where an orbit stays in the neighborhood, which we call 
the laminar lasting time distribution.
It is expected that by using the data-driven model 
the laminar lasting time distribution can be estimated in higher accuracy and 
in lower computational costs than by using the direct numerical simulation of the Navier--Stokes equation.
Here we study a macroscopic quantity of chaotic fluid flow in three dimensions under periodic boundary conditions~\cite{nakai_2018,nakai_2020}.

     \begin{figure}
		\subfigure[\hspace{-2mm}]{ \includegraphics[width=0.485\columnwidth,height=0.5\columnwidth]{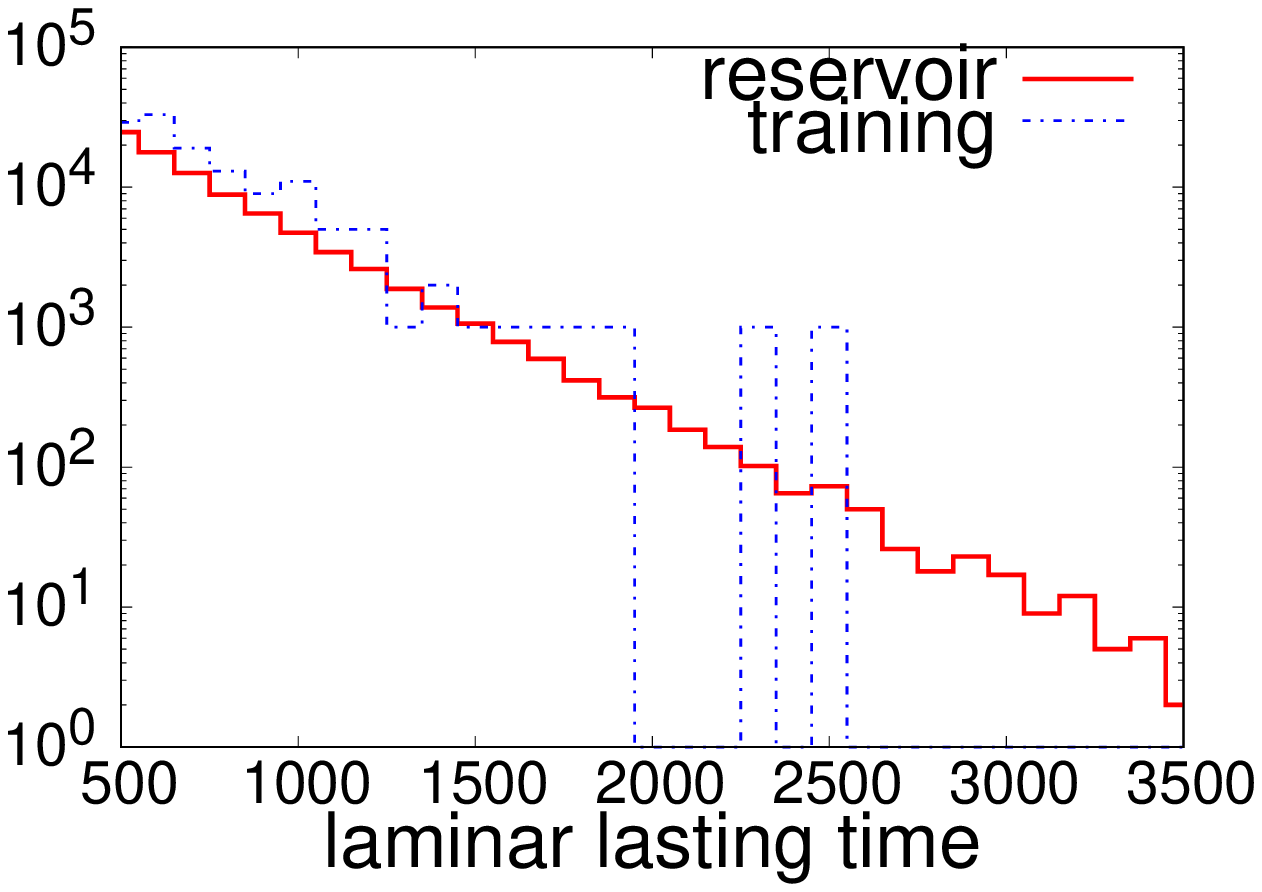}}%
		\subfigure[\hspace{-4mm}]{\includegraphics[width=0.485\columnwidth,height=0.5\columnwidth]{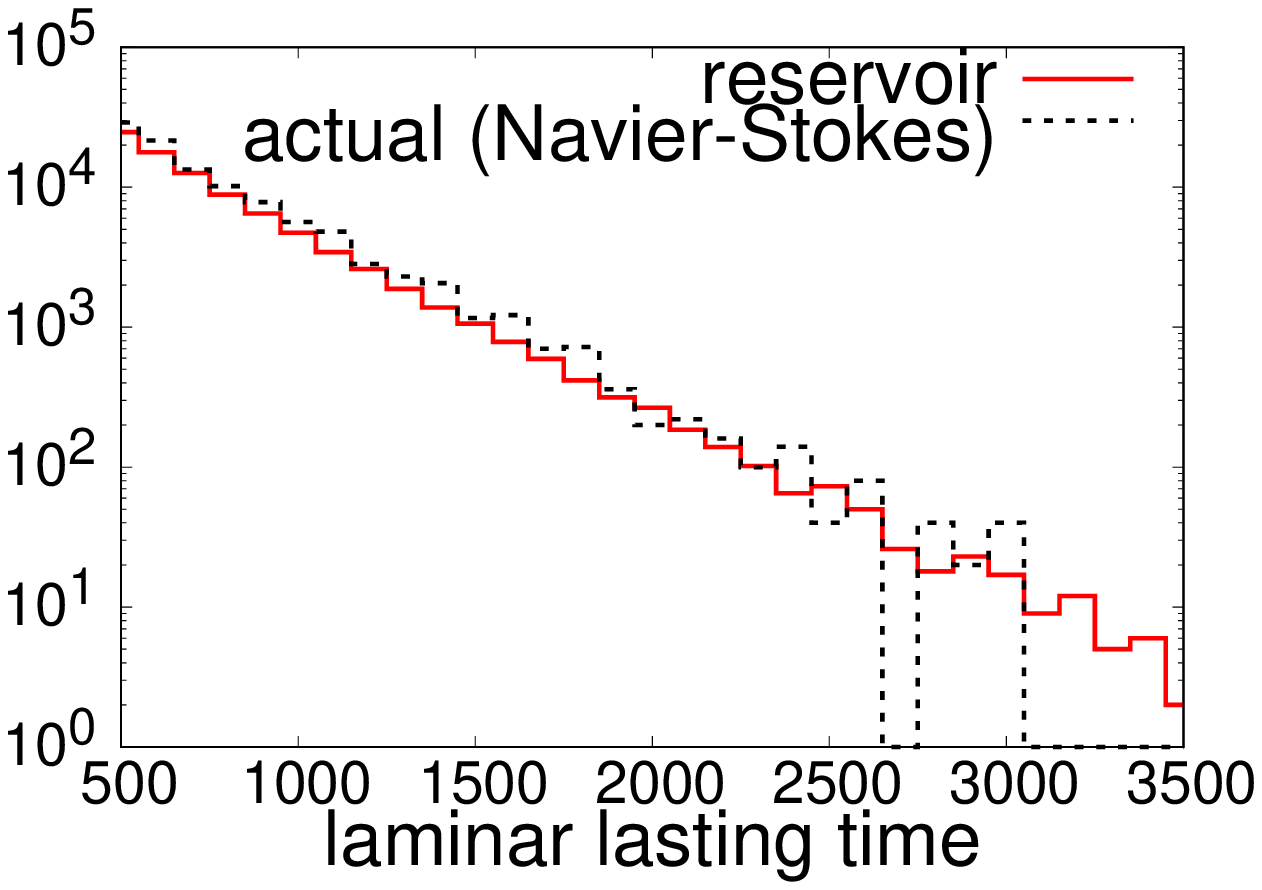}}
    \caption{\label{fig:lamlasttime}
    {\bf Laminar lasting time distribution of a fluid flow.}
      The laminar lasting time normalized distribution of a certain energy function $E(t)$ 
      of a fluid variable 
       (corresponding to $\tilde{E}(3,t)$ in \cite{nakai_2018}) estimated from the  trajectory of the data-driven model $(T_0 = 2 \times 10^8)$ is shown together with that from the short time training time series data $(T=0.001 \times T_0)$ (a) and with that by a long time actual time series data $(T=0.05 \times T_0)$ (b). The training data and actual data are calculated from a direct numerical simulation of the Navier--Stokes  equation. 
      $E(t)$ is the normalized variable (average 0, standard deviation 1) and 
      we consider the state is laminar when $|E(t)|<1.8$. 
}
	\end{figure}
The distribution shown in Fig.~\ref{fig:lamlasttime}
is generated from the very long trajectory of the data-driven model constructed 
by reservoir computing with a relatively low computational cost.
The detailed macroscopic dynamical structures 
can be determined using the data-driven model 
constructed from time series data without referring to microscopic behaviors. 
We hardly obtain these structures from a direct numerical simulations of the Navier--Stokes equation because of the high computational cost.
See the discrepancy in the distributions in  Fig.~\ref{fig:lamlasttime} (b).

It takes roughly $1/400$ of time to obtain a time-series of the energy functions $E(k)$ with the same time-lengths, when we use the model constructed by the reservoir computation. The Navier--Stokes equation is calculated by $13718$ dimensional ODEs with the $4$-stage Runge--Kutta method (time step $0.05$), whereas the model is calculated by 5000 dimensional map whose iterate corresponds to the time step $2$.

\section{VI. Concluding remarks.}
We have clarified by employing the time-series of the Lorenz 
system 
that a data-driven model using reservoir computing has quite similar dynamical system properties to those of the original Lorenz 
system, 
such as fixed points and their eigenvalues, periodic orbits, Lyapunov exponents and Lyapunov vectors.
It should be remarked that the fixed points extist far away from the training time-series data, but the corresponding points are found to exist nearby the original ones in the data-driven model.
We have also shown that 
the negative Lyapunov exponent computed not in the space of the reservoir state vector but in the space of output variables, and the 
degree of hyperbolicity measured by the angle between stable and unstable manifolds are shown to be quite similar 
to those of the original system. 
Qualitatively the same results are obtained 
for the R\"ossler system. 

 For a chaotic fluid flow we computed 
 the lasting time distribution staying in the neighborhood of a certain chaotic saddle showing laminar behavior by using the data-driven model. The model is constructed from a relatively short time-series data created from the direct numerical simulation of the Navier--Stokes equation.
 The obtained distribution cannot be computed from the direct computation of the Navier--Stokes equation because of its high computational cost.
 This result implies that a chaotic saddle can be reconstructed by the data-driven model.
\\
\indent{\bf Acknowledgements.} 
YS was supported by the JSPS KAKENHI Grant No.19KK0067 and No.21K18584.
KN was supported by the Project of President Discretionary Budget of TUMST. 
Part of the computation was supported by JHPCN (jh200020, jh210027), HPCI (hp200104, hp210072), and the Collaborative Research Program for Young $\cdot$ Women Scientists of ACCMS and IIMC, Kyoto University.

\bibliographystyle{apsrev4-2}

%

\end{document}